\documentstyle[amssymb,11pt]{article}
\begin{document}

\begin{center}{\bf  On The splitting of the Einstein field equations with respect to a general (1 + 3) threading of spacetime} \end{center}
\begin{center}{AUREL BEJANCU \\ DEPARTMENT OF MATHEMATICS \\ KUWAIT UNIVERSITY\\ P.O.Box 5969, Safat 13060\\ KUWAIT\\
 E-mail:aurel.bejancu@ku.edu.kw \vspace{5mm}\\
 HANI REDA FARRAN\\ DEPARTMENT OF MATHEMATICS \\ KUWAIT UNIVERSITY\\ P.O.Box 5969, Safat 13060\\ KUWAIT\\
 E-mail:hani.farran@ku.edu.kw} \end{center}

\begin{abstract}
{Based on a general (1+3) threading of the spacetime $(M, g)$, we obtain a new and simple splitting of a both the
Einstein field equations (EFE) and the conservation laws in $(M, g)$. As an application we obtain the splitting
of (EFE) in an almost FLRW universe with energy-momentum tensor of a perfect fluid. In particular, we state
the perturbation Friedman equations in an almost FLRW universe.} \end{abstract}

\medskip
\newpage
\section{Introduction}
The present paper is a continuation of the paper \cite{bc}, wherein an new approach on the (1+3) threading
of spacetime with respect to an arbitrary timelike vector field has been developed. The study in \cite{bc} refers
to Lorentz metrics given by (2.2) and subject to the condition that $\Phi$ is independent of the time coordinate.
We remove this condition, and therefore the results are valid for any Lorentz metric of a spacetime. Another
 important issue of the present paper is that the whole study is developed in the general setting of a spacetime
with a spatial distribution that is not necessarily integrable. The threading frames and coframes, the spatial
 tensor fields and the Riemannian spatial connection, are the main tools used throughout the paper. These geometric
 objects enable us to obtain a new and simple splitting of (EFE) and to apply it to the structure of an almost
FLRW universe. The new approach developed in this paper can be extended to a threading of higher-dimensional
 universes. In this respect we mention the paper \cite{b}, wherein a $(1+1+3)$-threading of a $5D$ universe has
been developed.\par
Now, we outline the content of the paper. In Section 2 we present the main geometric objects which constitute
the foundation of a general $(1+3)$ threading formalism of a spacetime $(M, g)$ with respect to an arbitrary
timelike vector field. We close this section with local expressions of the Levi-Civita connection $\nabla$ in terms
of spatial tensor field and of the local coefficients of the Riemannian connection $\bar{\nabla}$ (cf.(2.18)).
In Section 3 we state, for the first time in the literature, the spatial Bianchi identities in the general case
when the spatial distribution is not necessarily integrable (cf. (3.13), (3.14), (3.15)).  The structure equations
on $(M, g)$ induced by the $(1+3)$  threading formalism are presented in Section 4 (cf.(4.5),(4.6)). They play
an important role in the next sections, wherein we relate tensor fields on $M$ with spatial tensor fields. In
Section 5 we obtain simple expressions for the local components of the Ricci tensor of $(M, g)$ with respect to
the threading frame field , and for the scalar curvature (cf. (5.6), (5.7),
(5.8a), (5.11)). The splittings of both the Einstein gravitational tensor field and the energy-momentum tensor
field are stated in Section 6 (cf. (6.3), (6.11), (6.14), (6.16)). In Section 7 we obtain the spatial, mixed and
 temporal (EFE) (cf. (7.2), (7.3), (7.4)). Also, we state the equation (7.6) which in the particular case when
the treading is taken with respect to a unit vector field, becomes the well known Raychaudhuri-Ehlers equation. A
new splitting of conservation laws with respect to a general $(1+3)$ threading of spacetime is given in Section 8
(cf. (8.5), (8.6)). Also, we compare our results with what is known in the literature on this matter. Finally,
in Section 9 we apply the general theory developed in the paper to the $(1+3)$ threading of an almost FLRW
universe. We close the paper with conclusions and three appendixes.\newpage

\section{The $(1+3)$ threading formalism with respect to a non-normalized timelike vector field}

Recently, a new approach on the $(1+3)$ threading of spacetime with respect to a non-normalized timelike vector
field has been developed (cf. \cite{bc}). In the most general setting, that we explain in this section, we recall the
main geometric objects introduced in \cite{bc}. Also, we introduce the extrinsic curvature tensor field for the
 spatial distribution, and use it in the expressions of the Levi-Civita connection on a $4D$ spacetime.\par
Let $(M, g)$ be a $4D$ spacetime, and $\xi$ be a timelike vector field on $M$ that is not necessarily normalized.
Then, the  tangent  bundle $TM$ of $M$ admits the decomposition

$$ TM = VM\oplus SM.\eqno(2.1)$$
where $VM$ is the  {\it temporal distribution} spanned by $\xi$, and $SM$ is the {\it spatial distribution} that
is complementary orthogonal to $VM$ in $TM$.\par
Throughout the paper we use the ranges of indices: $i, j, k, ... \in\{1,2,3\}$ and $a, b, c, ...
\in\{0, 1, 2, 3\}$. Also, for any vector bundle $E$ over $M$ denote by $\Gamma(E)$ the ${\cal{F}}(M)$-module of
smooth sections of $E$, where ${\cal{F}}(M)$ is the algebra of smooth functions on $M$.\par
The congruence of curves that is tangent to $VM$, determines a  coordinate system $(x^a)$ on $M$ such that
$\xi = {\partial}/{\partial x^0}$. Next, we put

$$\begin{array}{lc}(a) \ \ \ \ \xi_0 = g\left(\frac{\partial }{\partial x^0},
\frac{\partial }{\partial x^0} \right) =
-\Phi^2,\ \ \ (b) \ \ \ \xi_i = g\left(\frac{\partial }{\partial x^i}, \frac{\partial }{\partial x^0}\right),
\ \vspace{2mm} \\  (c) \ \ \ \ \ \ \ g_{ij} = g\left(\frac{\partial }{\partial x^i},
\frac{\partial }{\partial x^j} \right), \end{array}\eqno(2.2)$$

wherein $\Phi$ is a non-zero function on $M$.\par
{\bf Remark 2.1} Note that in \cite{bc}, $\Phi$ was supposed to be independent of $x^0$. Here, we remove this
 condition on $\Phi$, and thus the results stated in the present paper are valid for any Lorentz metric $g$ on
$M$.\par
In this approach we use the {\it threading frame field} $\{\partial/\partial x^0, \delta/\delta x^i\}$ and the
{\it threading coframe field} $\{\delta x^0, dx^i\}$ defined as follows:

$$(a) \ \ \ \frac{\delta }{\delta x^i} = \frac{\partial }{\partial x^i} -
A_i\frac{\partial }{\partial x^0}, \ \ \ (b) \ \ \ \delta x^0 = dx^0 + A_idx^i,\eqno(2.3)$$
where we put
$$A_i = -\Phi^{-2}\xi_i.\eqno(2.4)$$
The Lie brackets of the vector fields from the threading frame are given by

$$(a) \ \ \left [\frac{\delta }{\delta x^j}, \frac{\delta }{\delta x^i}\right] =
2\omega_{ij}\frac{\partial }{\partial x^0},
\ \ (b) \ \ \left [\frac{\partial }{\partial x^0}, \frac{\delta }{\delta x^i}\right] =
a_i\frac{\partial }{\partial x^0},\eqno(2.5)$$
where we set

$$\begin{array}{lc}
(a) \ \ \ \omega_{ij} = \frac{1}{2}\left\{\frac{\delta A_j}{\delta x^i} - \frac{\delta A_i}{\delta x^j}\right\}
= \Phi^{-2}\left\{c_i\xi_j - c_j\xi_i + \frac{1}{2}\left(\frac{\delta\xi_i}{\delta x^j} - \frac{\delta\xi_j}{\delta x^i}
\right)\right\}, \vspace{2mm} \\ (b) \ \ \ \ c_i = \Phi^{-1}\frac{\delta\Phi}{\delta x^i}\ \ \ (c) \ \ \
a_i = -\frac{\partial A_i}{\partial x^0} = \Phi^{-2}\left\{\frac{\partial\xi_i}{\partial x^0} - 2\Psi\xi_i\right\},
 \vspace{2mm} \\ (d) \ \ \ \Psi = \Phi^{-1}\frac{\partial \Phi}{\partial x^0}. \end{array}\eqno(2.6)$$
Taking into account that the Levi-Civita connection $\nabla$ on $(M, g)$ is torsion-free, from (2.5a) we deduce that

$$\omega_{ij} = \frac{1}{2}\Phi^{-2}g\left(\nabla_{\frac{\delta}{\delta x^i}}\frac{\delta}{\delta x^i} -
\nabla_{\frac{\delta}{\delta x^j}}\frac{\delta}{\delta x^i}, \frac{\partial}{\partial x^0} \right).\eqno(2.7)$$
Thus, $\omega_{ij}, \ i,j \in \{1,2,3\}$, define the {\it vorticity tensor field} on $(M, g)$. By using the
Jacobi identity

$$\left [\left [X, Y\right ], Z\right ] + \left [\left [Y, Z\right ], X\right ] + \left [\left [Z, X\right ], Y\right ]
= 0, \ \ \ \forall\ X, Y, Z \in \Gamma(TM),$$
we deduce that the vorticity tensor field satisfies the identities

$$(a) \frac{\partial \omega_{ij}}{\partial x^0} = \frac{1}{2}\left\{\frac{\delta a_i}{\delta x^j} -
\frac{\delta a_j}{\delta x^i} \right\}, \ \ \ (b) \ \ \ \sum_{(i,j,k)}\left\{\frac{\delta\omega_{ij}}{\delta x^k} -
\omega_{ij}a_k\right\} = 0, \eqno(2.8)$$
where $\displaystyle{\sum_{(i,j,k)}}$ is the cyclic sum with respect to $(i,j,k)$. \par
Now, we denote by $\bar{g}_{ij}$ the local components of the Riemannian metric $\bar{g}$ induced by $g$ on $SM$,
with respect ti the basis $\{\delta/\delta x^i\}$ in $\Gamma(SM)$, and obtain

$$\bar{g}_{ij} = \bar{g}(\frac{\delta }{\delta x^i}, \frac{\delta }{\delta x^j}) = g_{ij} + \Phi^2A_iA_j =
g_{ij} + \Phi^{-2}\xi_i\xi_j.\eqno(2.9)$$
 Then the lone element of $g$ is expressed in terms of threading coframe $\{\delta x^0, dx^i\}$ as follows:

$$ds^2 = - \Phi^2(\delta x^0)^2  + \bar{g}_{ij}dx^idx^j .\eqno(2.10) $$
Also, we define the {\it expansion tensor field} $\Theta_{ij}$, the {\it expansion function} $\Theta$ and the
{\it shear tensor field} $\sigma_{ij}$ as follows:

$$(a) \ \ \Theta_{ij} = \frac{1}{2}\frac{\partial \bar{g}_{ij}}{\partial x^0}, \ \ \ (b) \ \ \
 \Theta = \Theta_{ij}\bar{g}^{ij}, \ \ \ (c) \ \ \ \sigma_{ij} = \Theta_{ij} - \frac{1}{3}\Theta
\bar{g}_{ij}. \eqno(2.11)$$

Raising and lowering indices $i,j,k,...$ are performed by using $\bar{g}^{ij}$ and $\bar{g}_{ij}$, as in the
following examples:

$$\omega_j^k = \bar{g}^{ki}\omega_{ij}, \ \ \ \omega^{kh} = \bar{g}^{ki}\bar{g}^{hj}\omega_{ij}, \ \ \ \omega_{ij}
= \bar{g}^{ik}\omega_j^k, \ \ \ \omega_{ij} =  \bar{g}_{ik}\bar{g}_{jh}\omega^{kh}.$$
The expansion and vorticity tensor fields enable us to define the {\it extrinsic curvature tensor field} $K$
of the spatial distribution by its local components

$$K_j^h = \Theta_j^h + \Phi^2\omega_j^h,\eqno(2.12) $$
or equivalently by

$$K_{ij} = \Theta_{ij} + \Phi^2\omega_{ij}.\eqno(2.13) $$
By using (2.5a) and (2.13), we see that $K$ is a symmetric tensor field if and only if $SM$ is integrable.
\vspace{2mm}\par
{\bf Remark 2.2} The extrinsic curvature tensor field was intensively used in the $(3+1)$ decomposition of the
 spacetime (cf.\cite{m}, pp 509-516). As far as we know, the tensor field $K$ given by (2.12) or (2.13) is
considered here for the first time in a study of the $(1+3)$ threading of spacetime.\vspace{2mm}\par
Next, in order to justify the tensorial meaning of the above quantities, we define a {\it spatial tensor field} $T$ of
type $(p, q)$ on $M$, as an ${\cal{F}}(M)$-multilinear mapping

$$T:\Gamma(SM^{\star})^p\times \Gamma(SM)^q\longrightarrow {\cal{F}}(M),$$
where $SM^{\star}$ is the dual vector bundle to $SM$. The local components of $T$ with respect to a threading
frame and coframe, are given by

$$T_{i...}^{k...} = T(dx^k,...,\frac{\delta }{\delta x^i},...),$$
and satisfy

$$T_{i \cdots }^{k\cdots }\frac{\partial\tilde{x}^h}{\partial x^k} = \tilde{T}^{h \cdots }_{j\cdots }
 \frac{\partial\tilde{x}^j}{\partial x^i},$$
 with respect to the coordinate transformations $\tilde{x}^a = \tilde{x}^a(x^0, x^i)$ on $M$. As examples:
 $\{\omega_{ij}, \bar{g}_{ij}, \theta_{ij}, \sigma_{ij}, K_{ij} \}$ and $\{a_i, c_i\}$ define spatial tensor fields
of type $(0, 2)$ and $(0.1)$ respectively.\par
An important geometric object is the {\it Riemannian spatial connection}, which is a metric linear connection
 $\bar{\nabla}$ on the spatial distribution, given by

$$(a) \ \ \ \bar{\nabla}_{X}{\cal{S}}Y = {\cal{S}}{\nabla}_X{\cal{S}}Y, \ \ \ \forall \ \ X,Y
\in \Gamma(TM),\eqno(2.14)$$
where ${\cal{S}}$ is the projection morphism of $TM$ on $SM$ with respect to the decomposition (2.1).
 Locally, $ \bar{\nabla}$ is given by

$$(a) \ \ \bar{\nabla}_{\frac{\delta }{\delta x^j}}\frac{\delta }{\delta x^i} =
\bar{\Gamma}^{\; \;k}_{i\ \;j}\frac{\delta }{\delta x^k}, \ \ \ (b) \ \ \ \
\bar{\nabla}_{\frac{\partial }{\partial x^0}}\frac{\delta }{\delta x^i} =
K^k_i\frac{\delta }{\delta x^k}.\eqno(2.15) $$
where we put
$$\bar{\Gamma}^{\; \;k}_{i\ \;j} = \frac{1}{2}\bar{g}^{kh}\left\{\frac{\delta \bar{g}_{hj}}{\delta x^i} +
\frac{\delta \bar{g}_{hi}}{\delta x^j} - \frac{\delta \bar{g}_{ij}}{\delta x^h}\right\}. \eqno(2.16) $$
Throughout the paper, the covariant derivatives defined by $\bar{\nabla}$ will be denote by a vertical bar
$"\vert"$. As an example, for a spatial tensor field  $T = (T^i_j)$ we have

$$ \begin{array}{c}(a) \ \ \
T_{i\vert_{k}}^{j} = \frac{\delta T_i^j}{\delta x^k} + T_i^h\bar{\Gamma}^{\; \;j}_{h\ \;k}
  - T_{h}^{j}\bar{\Gamma}^{\; \;h}_{i\ \;k},\vspace{3mm}\\ (b) \ \ \
 T_{i\vert_{0}}^j = \frac{\partial T_i^j}{\partial x^0} + T_i^hK^j_h - T_h^jK^h_i.\end{array}  \eqno(2.17)$$
A covariant derivative as in (2.17a) (resp. (2.17b)) is called a {\it spatial covariant derivative} (resp. {\it
temporal covariant derivative}) of the spatial tensor field $T$.\par
Finally, by direct calculations, using the Riemannian spatial connection and the above spatial tensor fields,
we express the Levi-Civita connection $\nabla$ on $(M, g)$ as follows:

$$\begin{array}{lc}
(a) \ \ \nabla_{\frac{\delta }{\delta x^j}}\frac{\delta }{\delta x^i} =  \bar{\Gamma}^{\; \;k}_{i\ \;j}
\frac{\delta }{\delta x^k} + \left(\omega_{ij} + \Phi^{-2}\Theta_{ij}\right)\frac{\partial }{\partial x^0}
\vspace{3mm}\\ \hspace*{22mm} =  \bar{\Gamma}^{\; \;k}_{i\ \;j}\frac{\delta }{\delta x^k}
+ \Phi^{-2}K_{ij}\frac{\partial }{\partial x^0}, \vspace{4mm}\\
(b) \ \ \nabla_{\frac{\partial }{\partial x^0}}\frac{\delta }{\delta x^i} = \left(\Theta_i^k
+ \Phi^2\omega_i^k\right)\frac{\delta }{\delta x^k} + b_i\frac{\partial }{\partial x^0}
\vspace{3mm}\\ \hspace*{22mm} =  K^k_i\frac{\delta }{\delta x^k}
+ b_i\frac{\partial }{\partial x^0}, \vspace{4mm}\\
(c) \ \ \nabla_{\frac{\delta }{\delta x^i}}\frac{\partial }{\partial x^0} =
\left(\Theta_i^k + \Phi^2\omega_i^k\right)\frac{\delta }{\delta x^k} + c_i\frac{\partial }{\partial x^0}
\vspace{3mm}\\ \hspace*{22mm} =  K^k_i\frac{\delta }{\delta x^k}
+ c_i\frac{\partial }{\partial x^0}, \vspace{4mm}\\
(d) \ \ \nabla_{\frac{\partial }{\partial x^0}}\frac{\partial }{\partial x^0} =
 \Phi^2b^k\frac{\delta }{\delta x^k} + \Psi\frac{\partial }{\partial x^0}, \end{array}\eqno(2.18)$$
where we put

$$b_i = a_i + c_i, \ \ \ \ i \in \{1, 2, 3\}.\eqno(2.19) $$
{\bf Remark 2.3} It is worth mentioning that all the equations we state in the paper are expressed in
terms of spatial tensor fields and their covariant derivatives defined by the Riemannian connection.\par
 {\bf Remark 2.4} As the $(1 + 3)$ threading of spacetime considered in this paper contains as a particular case
the $(1 + 3)$ threading with respect to a unit timelike vector field, we call it the {\it general $(1 + 3)$
threading of spacetime}. The advantage of this general setting on the splitting of spacetime is that it can be
applied to any Lorentz metric of a spacetime.

 \section{ Bianchi identities for the Riemannian spatial connection}
In earlier Literature on the $(1 + 3)$ threading of spacetime  we find the so called three-dimensional
derivative operator (cf.(4.19) of \cite{ej}). With respect to this operator we have the following remarks:\par
(i) It is neither a linear connection on $M$, nor a linear connection on $SM$.\par
(ii) As a consequence of (i), for the general case when $SM$ is not integrable, then a curvature tensor field for
this operator could be not defined.\par
Contrary to this situation, $\bar{\nabla}$ given by (2.14) is a metric linear connection on the vector bundle
$SM$, and therefore it has a curvature tensor field $\bar{R}$ given by

$$\begin{array}{lc}
\bar{R}(X, Y, Z) = \bar{\nabla}_X\bar{\nabla}_Y{\cal{S}}Z - \bar{\nabla}_Y\bar{\nabla}_X{\cal{S}}Z
- \bar{\nabla}_{[X, Y]}{\cal{S}}Z,\vspace{2mm} \\ \hfill \forall \ X, Y, Z \in \Gamma(TM).\end{array}\eqno(3.1)$$
Locally, we put

$$\begin{array}{l}(a) \ \ \ \bar{R}(\frac{\delta }{\delta x^k}, \frac{\delta }{\delta x^j},
\frac{\delta }{\delta x^i}) = \bar{R}_{i\ jk}^{\; h}\frac{\delta }{\delta x^h}, \vspace{3mm}\\
 (b) \ \ \ \bar{R}(\frac{\delta }{\delta x^k}, \frac{\partial }{\partial x^0}, \frac{\delta }{\delta x^i}) =
\bar{R}_{i\  0 k}^{\;h}\frac{\delta }{\delta x^h}, \end{array}\eqno(3.2)$$
and by using (3.1), (3.2), (2.15), and (2.5), we obtain

$$\begin{array}{lc}(a) \ \ \ \bar{R}^{\; h}_{i\ jk} = \frac{\delta\bar{\Gamma}_{i\;\;j}^{\; h}}{\delta x^k} -
\frac{\delta\bar{\Gamma}_{i\;\;k}^{\; h}}{\delta x^j} + \bar{\Gamma}_{i\;\;j}^{\;l}\bar{\Gamma}_{l\;\;k}^{\;h}
- \bar{\Gamma}_{i\;\;k}^{\; l}\bar{\Gamma}_{l\;\;j}^{\; h} \vspace{3mm}\\ \hspace*{12mm}
- 2K_i^h\omega_{jk}, \vspace{3mm}\\
(b) \ \ \ \bar{R}^{\; h}_{i\ 0k} = K_{i\ \vert k}^{\;h} - \frac{\partial \bar{\Gamma}_{i\;\;k}^{\;h}}
{\partial x^0} + K_i^ha_k.\end{array}\eqno(3.3)$$
Sine $\bar{\nabla}$ is a metric linear connection, we have

$$\bar{R}_{il0k} + \bar{R}_{li0k} = 0, \eqno(3.4)$$
where we put

$$\bar{R}_{il0k} = \bar{g}_{lh}\bar{R}^{\; h}_{i\ 0k}. \eqno(3.5) $$
As a consequence of (3.4) we deuce that $\bar{R}^{\; i}_{i\ 0k} = 0,$ which implies

$$\Theta_{\vert k} = \frac{\partial \Theta}{\partial x^k} - A_k\frac{\partial \Theta}{\partial x^0} =
\frac{\bar{\Gamma}_{i\;\;k}^{\; i}}{\partial x^0} - \Theta a_k, \eqno(3.6)$$
via (3.3b) and (2.12). Thus, {\it in any cosmological model of a $4D$ universe the expansion function must satisfy
the system of $PDE$ given by} (3.6).\par
{\bf Remark 3.1} Note that $\bar{R}^{\; h}_{i\ jk}$ and $\bar{R}^{\; h}_{i\ 0k}$ define spatial tensor fields of type
$(1, 3)$ and $(1, 2)$, respectively. Also, from (3.3b) we see that $\partial{\bar{\Gamma}_{i\;\;k}^{\; h}}/
{\partial x^0}$ define a spatial tensor field of type $(1, 2)$. However, $\bar{\Gamma}_{i\;\;k}^{\; h}$ do not define
a spatial tensor field.\par
{\bf Remark 3.2} Comparing (3.3a) with (15.4) from \cite{a}, we see that the so called Zelmanov curvature tensor
field is given by the first four terms from (3.3a). Moreover, from (15.5) of \cite{a} we see that such a tensor
field becomes a curvature tensor field, if and only if, $SM$ is an integrable distribution.\par
Next, we extend the Riemannian spatial connection $\bar{\nabla}$ on $SM$ to a linear connection $\tilde{\nabla}$
on $M$ given by

$$\tilde{\nabla}_XY =  \bar{\nabla}_X{\cal{S}}Y + {\cal{T}}\nabla_X{\cal{T}}Y, \ \ \ \forall \ X, Y
\ \in\Gamma(TM),\eqno(3.7)$$
where ${\cal{T}}$ is the projection morphism of $TM$ on $VM$ with respect to (2.1). Clearly, $\tilde{\nabla}$
 coincides with $\bar{\nabla}$ on $SM$ and therefore locally is given by (2.15) and

$$(a) \ \ \ \tilde{\nabla}_{\frac{\delta }{\delta x^i}}\frac{\partial }{\partial x^0} = c_i\frac{\partial }
{\partial x^0}, \ \ \ (b) \ \ \ \tilde{\nabla}_{\frac{\delta }{\delta x^0}}\frac{\partial }{\partial x^0} =
\Psi\frac{\partial }{\partial x^0}.\eqno(3.8)$$
We recall that the torsion and curvature tensor fields of $\tilde{\nabla}$ are given by

$$\begin{array}{lc}(a) \ \ \ \tilde{T}(X, Y) = \tilde{\nabla}_XY -  \tilde{\nabla}_YX - [X, Y],\vspace{2mm}\\
(b) \ \ \ \tilde{R}(X, Y, Z) = \tilde{\nabla}_X\tilde{\nabla}_YZ - \tilde{\nabla}_Y\tilde{\nabla}_XZ
- \tilde{\nabla}_{[X, Y]}Z.\end{array}\eqno(3.9)$$
Then, by direct calculations, using (3.9), (2.15), (3.8), (2.5), (2.18c) and (3.2), we deduce that

$$\begin{array}{l}(a) \ \ \ \tilde{T}(\frac{\delta }{\delta x^j}, \frac{\delta }{\delta x^i}) = -2\omega_{ij}
\frac{\partial }{\partial x^0},\ \ \ \vspace{3mm}\\(b) \ \ \ \tilde{T}(\frac{\partial }{\partial x^0},
\frac{\delta }{\delta x^i})
 = K_i^j\frac{\delta }{\delta x^j} - b_i\frac{\partial }{\partial x^0}, \vspace{3mm}\\
 (c) \ \ \ \tilde{R}(\frac{\delta }{\delta x^h}, \frac{\delta }{\delta x^k}, \frac{\delta }{\delta x^i}) =
\bar{R}_{i\ kh}^{\;j}\frac{\delta }{\delta x^j},\ \ \ \vspace{3mm}\\
 (d) \ \ \ \tilde{R}(\frac{\delta }{\delta x^k}, \frac{\partial }{\partial x^0}, \frac{\delta }{\delta x^i}) =
\bar{R}_{i\  0 k}^{\;j}\frac{\delta }{\delta x^j}. \end{array}\eqno(3.10)$$
Now, in order to find some Bianchi identities for the Riemannian spatial connection, we recall that the
Bianchi identities for the linear connection $\tilde{\nabla}$ are given by (cf.\cite{s}, p.135)

$$\begin{array}{lc}(a) \ \ \ \displaystyle{\sum_{(X,Y,Z)}}\{(\tilde{\nabla}_X\tilde{T})(Y, Z)
+ \tilde{T}(\tilde{T}(X, Y), Z)
\vspace{2mm}\\ \hspace*{22mm}- \tilde{R}(X,Y, Z) \} = 0, \vspace{2mm}\\
(b) \ \ \ \displaystyle{\sum_{(X,Y,Z)}}\{(\tilde{\nabla}_X\tilde{R})(Y, Z, U) +
\tilde{R}(\tilde{T}(X, Y), Z, U) \} = 0,\end{array}\eqno(3.11)$$
where $\displaystyle{\sum_{(X,Y,Z)}}$ is the cyclic sum with respect to $(X,Y,Z)$.
In order to use (3.11a), we note that

$$\begin{array}{lc} (a) \ \ \
\tilde{T}(\tilde{T}(\frac{\delta }{\delta x^k}, \frac{\delta }{\delta x^j}), \frac{\delta }{\delta x^i}) =
-2K_i^h\omega_{jk}\frac{\delta }{\delta x^h} + 2b_i\omega_{jk}\frac{\partial }{\partial x^0}, \vspace{2mm}\\
(b) \ \ \ (\tilde{\nabla}_{\frac{\delta }{\delta x^k}}\tilde{T})(\frac{\delta }{\delta x^j},
\frac{\delta }{\delta x^i}) = -2(\omega_{ij\vert k} + \omega_{ij}c_k)\frac{\partial }{\partial x^0}.
\end{array}\eqno(3.12)$$
Then, take $X = \delta/\delta x^k$, $Y = \delta/\delta x^j$, $Z = \delta/\delta x^i$ in (3.11a) and by using
(3.12) and (3.10a), we infer that the spatial component in (3.11a) is expressed as follows:

$$ \sum_{(i,j,k)}\{\bar{R}^{\;\;h}_{i\ jk} + 2K_i^h\omega_{jk}\} = 0.\eqno(3.13)$$
Taking temporal part in (3.11a) and any other triplet $(X, Y, Z)$ from the threading frame we obtain the
identities from (2.8)\par
Next, take $X = \delta/\delta x^k$, $Y = \delta/\delta x^j$, $Z = \delta/\delta x^i$ and $U = \delta/\delta x^h$
in (3.11b), and by using (3.7), (3.10a), (3.10c) and (3.10d), we obtain

$$ \sum_{(i,j,k)}\{\bar{R}^{\;\;l}_{h\ ij\vert k} + \bar{R}^{\;\;l}_{h\ 0i\vert k}\omega_{jk}\} = 0.\eqno(3.14)$$
Finally, take $X = \partial/\partial x^0$, $Y = \delta/\delta x^j$, $Z = \delta/\delta x^i$ and $U = \delta/\delta x^h$
in (3.11b), and by using (3.7), (3.10) and (2.18c), we deduce the identity

$$\bar{R}^{\;\;l}_{h\ ij\vert 0} + \bar{R}^{\;\;l}_{h\ 0i\vert j} - \bar{R}^{\;\;l}_{h\ 0j\vert i} +
\bar{R}^{\;\;l}_{h\ ik}K_j^k - \bar{R}^{\;\;l}_{h\ jk}K_i^k + a_j\bar{R}^{\;\;l}_{h\ 0i} - a_i\bar{R}^{\;\;l}_{h\ 0j}
= 0.\eqno(3.15)$$
The other identities obtained from (3.11b) are either trivial or they do not involve the curvature tensor $\bar{R}$ of
$\bar{\nabla}$. Thus, we are entitled to call (3.13), (3.14) and (3.15) the {\it Bianchi identities} for
the Riemannian spatial connection.\par
We close the section with some comments on these identities. As for as we know, the above Bianchi identities are
stated here for the first time in the literature. They represent a generalization of usual Bianchi identities on
a 3-dimensional Riemannian manifold. Indeed, if the spatial distribution is integrable, that is, the vorticity
tensor field vanishes identically on $M$, then (3.13) and (3.14) become

$$ \sum_{(i,j,k)}\{\bar{R}^{\;\;l}_{i\ jk}\} = 0, \ \ \ \ \sum_{(i,j,k)}\{\bar{R}^{\;\;l}_{h\ ij\vert k}\} = 0.$$
which are the well-known Bianchi identities on the 3-dimensional leaves of $SM$. Moreover, in this case, by
using (2.12) and (3.3) we deduce that

$$\begin{array}{lc}(a) \ \ \ K_i^j = \Theta_i^j, \ \ \ (b) \ \ \ \bar{R}^{\; j}_{i\ 0h} = \Theta_{i\vert h}^{j}
- \frac{\partial \bar{\Gamma}_{i\;\;h}^{\;j}}{\partial x^0} + \Theta_i^ja_h,\vspace{3mm}\\
(c) \ \ \ \bar{R}^{\; j}_{i\ kh} = \frac{\delta\bar{\Gamma}_{i\;\;k}^{\; j}}{\delta x^h}
- \frac{\delta\bar{\Gamma}_{i\;\;h}^{\; j}}{\delta x^k} +
\bar{\Gamma}_{i\;\;k}^{\;l}\bar{\Gamma}_{l\;\;h}^{\;j}
- \bar{\Gamma}_{i\;\;h}^{\; l}\bar{\Gamma}_{l\;\;k}^{\; j}. \end{array}\eqno(3.16)$$
Finally, the identity (3.15) becomes

$$\bar{R}^{\;\;l}_{h\ ij\vert 0} + \bar{R}^{\;\;l}_{h\ 0i\vert j} - \bar{R}^{\;\;l}_{h\ 0j\vert i} +
\bar{R}^{\;\;l}_{h\ ik}\Theta_j^k - \bar{R}^{\;\;l}_{h\ jk}\Theta_i^k + a_j\bar{R}^{\;\;l}_{h\ 0i}
- a_i\bar{R}^{\;\;l}_{h\ 0j} = 0.\eqno(3.17)$$

\section{Structure equations induced by the $(1 + 3)$ threading of spacetime}

Let $(M, g)$ be a $4D$ spacetime and $\nabla$ be the Levi-Civita connection defined by the Lorentz metric $g$. Denote by
$R$ the curvature tensor field of $\nabla$ given by (3.9b) wherein we remove the tilde. Then, consider the following
local components of $R$ with respect to the threading frame $\{\partial/\partial x^0, \delta/\delta x^i\}:$

$$\begin{array}{l}(a) \ \ \ R(\frac{\delta }{\delta x^k}, \frac{\delta }{\delta x^j}, \frac{\delta }{\delta x^i}) =
R_{i\ jk}^{\;h}\frac{\delta }{\delta x^h} + R_{i\ jk}^{\;0}\frac{\partial }{\partial x^0},\ \ \ \vspace{3mm}\\
(b) \ \ \ R(\frac{\delta }{\delta x^k}, \frac{\partial }{\partial x^0},
\frac{\delta }{\delta x^i}) =
R_{i\  0 k}^{\;h}\frac{\delta }{\delta x^h} + R_{i\ 0k}^{\;0}\frac{\partial }{\partial x^0}. \end{array}\eqno(4.1)$$
Now, comparing $(A_4)$ and $(A_8)$ from Appendix $A$ with (4.1a) and (4.1b) respectively, we obtain

$$\begin{array}{l}(a) \ \ \
R_{i\ jk}^{\;\;h} = \bar{R}_{i\ jk}^{\;\;h} + \Phi^{-2}\left(K_{ij}K_k^h - K_{ik}K_j^h \right), \vspace{2mm}\\
(b) \ \ \ R_{i\ jk}^{\;\;0} = \Phi^{-2}\left(K_{ij\vert k} - K_{ik\vert j}+ K_{ik}c_j -K_{ij}c_k\right) -
2b_i\omega_{jk}, \vspace{2mm}\\
(c) \ \ \ R_{i\ 0k}^{\;\;h} = \bar{R}_{i\ 0k}^{\;\;h} + b_iK_k^h - b^hK_{ik}, \vspace{2mm}\\
(d) \ \ \  R_{i\ 0k}^{\;\;0} = b_{i\vert k} + b_ib_k - \Phi^{-2}\left(K_{ik\vert 0} + K_{ij}K^j_k -
\Psi K_{ik}\right).\end{array} \eqno(4.2)$$
By using (2.13), (2.6) and (2.6d), we deduce that

$$\begin{array}{c}(a) \ \ \ K_{ij\vert k} = \Theta_{ij\vert k} + \Phi^2\left(\omega_{ij\vert k} +
2\omega_{ij}c_k \right), \vspace{2mm} \\ (b) \ \
K_{ij\vert 0} = \Theta_{ij\vert 0} + \Phi^2\left(\omega_{ij\vert 0} +
2\Psi\omega_{ij}\right).\end{array}\eqno(4.3)$$
Taking account of (4.3), (2.12) and (2.13) into (4.2), we infer that

$$\begin{array}{l}(a) \ \ \
R_{i\ jk}^{\;\;h} = \bar{R}_{i\ jk}^{\;\;h} + \left(\omega_{ij} + \Phi^{-2}\Theta_{ij}\right)
\left(\Theta_k^h + \Phi^2\omega_k^h \right) \vspace{2mm}\\ \hspace*{22mm} - \left(\omega_{ik}
+ \Phi^{-2}\Theta_{ik}\right)\left(\Theta_j^h + \Phi^2\omega_j^h\right),\vspace{2mm}\\
(b) \ \ \ R_{i\ jk}^{\;\;0} = \Phi^{-2}\left(\Theta_{ij\vert k} - \Theta_{ik\vert j}+ \Theta_{ik}c_j -
\Theta_{ij}c_k\right) \vspace{2mm}\\ \hspace*{22mm} +\ \omega_{ij\vert k} - \omega_{ik\vert j}+ \omega_{ij}c_k -
\omega_{ik}c_j - 2b_i\omega_{jk}, \vspace{2mm}\\
(c) \ \ \ R_{i\ 0k}^{\;\;h} = \bar{R}_{i\ 0k}^{\;\;h} + b_i\Theta_k^h - b^h\Theta_{ik} +
\Phi^2\left(b_i\omega_k^h - b^h\omega_{ik} \right), \vspace{2mm}\\
(d) \ \ \  R_{i\ 0k}^{\;\;0} = b_{i\vert k} + b_ib_k - \Phi^{-2}\left(\Theta_{ik\vert 0} -
\Psi \Theta_{ik}\right)\vspace{2mm}\\ \hspace*{22mm} - \omega_{ik\vert 0} - \Psi\omega_{ik}- \left(\omega_{ij}
+ \Phi^{-2}\Theta_{ij}\right)\left(\Theta_k^j + \Phi^2\omega_k^j \right).\end{array} \eqno(4.4)$$
Next, by using the local components of the curvature tensor fields $R$ and $\bar{R}$ of type $(0, 4)$ (see $(A_9)$
and $(A_{10})$), from (4.2) and (4.4) we obtain

$$\begin{array}{l}(a) \ \ \
R_{iljk} = \bar{R}_{iljk} + \Phi^{-2}\left(K_{ij}K_{lk} - K_{ik}K_{lj} \right), \vspace{2mm}\\
(b) \ \ \ R_{i0jk} = K_{ik\vert j} - K_{ij\vert k}+ K_{ij}c_k -K_{ik}c_j +
2\Phi^2b_i\omega_{jk}, \vspace{2mm}\\
(c) \ \ \ R_{il0k} = \bar{R}_{il0k} + b_iK_{lk} - b_lK_{ik}, \vspace{2mm}\\
(d) \ \ \  R_{i00k} = K_{ik\vert 0} + K_{ij}K^j_k - \Psi K_{ik} - \Phi^{-2}\left(b_{i\vert k}
+ b_ib_k\right),\end{array} \eqno(4.5)$$
and

$$\begin{array}{l}(a) \ \ \
R_{iljk} = \bar{R}_{ijk} + \Phi^2\left\{\left(\omega_{ij} + \Phi^{-2}\Theta_{ij}\right)
\left(\omega_{lk} + \Phi^{-2}\Theta_{lk}\right)\right. \vspace{2mm}\\ \hspace*{22mm} - \left.\left(\omega_{ik}
+ \Phi^{-2}\Theta_{ik}\right)\left(\omega_{lj} + \Phi^{-2}\Theta_{lj}\right)\right\},\vspace{2mm}\\
(b) \ \ \ R_{i0jk} = \Theta_{ik\vert j} - \Theta_{ij\vert k}+ \Theta_{ij}c_k -
\Theta_{ik}c_j \vspace{2mm}\\ \hspace*{22mm} + \Phi^{2}\left\{\omega_{ik\vert j} - \omega_{ij\vert k}+ \omega_{ik}c_j -
\omega_{ij}c_k + 2b_i\omega_{jk}\right\}, \vspace{2mm}\\
(c) \ \ \ R_{il0k} = \bar{R}_{il0k} + b_i\Theta_{lk} - b_l\Theta_{ik} +
\Phi^2\left(b_i\omega_{lk} - b_l\omega_{ik} \right), \vspace{2mm}\\
(d) \ \ \  R_{i00k} = \Theta_{ik\vert 0} - \Psi\Theta_{ik} + \Phi^{2}\left\{\omega_{ik\vert 0} + \Psi\omega_{ik}
\right.\vspace{2mm}\\ \hspace*{22mm}\left. + \left(\omega_{ij} + \Phi^{-2}\Theta_{ij}\right)\left(\Theta_{k}^j
+ \Phi^{2}\omega_{k}^j\right) - b_{i\vert k} - b_ib_k\right\}.\end{array} \eqno(4.6)$$
With the theory of hypersurfaces of the spacetime in mind, we call (4.5a) and (4.6a) \{resp. (4.5b), (4.6b) and
(4.6c)\} the
{\it Gauss equations} (resp. {\it Codazzi equations}) for the spatial distribution $SM$ in the ambient space
$(M, g)$. Also, all the equations from both (4.5) and (4.6) will be called {\it structure equations} induced by the
$(1+3)$ threading formalism. They have an important role in the next sections.\par
Now, taking into account the symmetries of $R$ we deduce some identities for $\bar{R}$ and for kinematic
 quantities. First, using well known identities for $R$ and taking into account (4.5a) and (4.6a), we obtain
the following identities for $\bar{R}$:

$$\begin{array}{lc}
(a) \ \ \ \ \bar{R}_{iljk} + \bar{R}_{ilkj} = 0, \vspace{2mm} \\ (b) \ \ \ \ \bar{R}_{iljk} + \bar{R}_{lijk} = 0,
\vspace{2mm}\\ (c) \ \ \ \ \bar{R}_{iljk} - \bar{R}_{jkil} = \Phi^{-2}\left\{K_{ik}K_{lj}  + K_{ji}K_{kl} -
K_{ij}K_{lk} - K_{jl}K_{ki}\right\}\vspace{2mm}\\ \hspace*{22mm} = 2\left\{\Theta_{ik}\omega_{lj}  +
\Theta_{lj}\omega_{ik} + \Theta_{ij}\omega_{kl} + \Theta_{kl}\omega_{ji}\right\}.\end{array}\eqno(4.7)$$
Also, taking into account that

$$R_{i0jk} = -R_{jk0i},$$
and by using (4.5b), (4.5c), (4.6b) and (4.6c), we deduce that

$$\begin{array}{lc}
\bar{R}_{jk0i} = K_{ij\vert k} - K_{ik\vert j} + K_{ik}c_j - K_{ij}c_k
- 2\Phi^2b_i\omega_{jk} - b_jK_{ki} + b_kK_{ji}\vspace{2mm}\\ \hspace*{10mm} = \
\Theta_{ij\vert k} - \Theta_{ik\vert j} + \Theta_{ij}a_k - \Theta_{ik}a_j \vspace{2mm}\\ \hspace*{12mm}
+ \ \Phi^{2}\left\{\omega_{ij\vert k} - \omega_{ik\vert j} + \omega_{ik}a_j -
\omega_{ij}a_k - 2b_i\omega_{jk}\right\}.\end{array}\eqno(4.8)$$
Finally, using the identity

$$R_{i00k} = R_{k00i},$$
and taking the symmetric and skew-symmetric parts in(4.5d) and (4.6d), we infer that

$$\begin{array}{lc}(a) \ \ \  R_{i00k} = \frac{1}{2}\left\{K_{ik\vert 0} + K_{ki\vert 0} + K_{ij}K^j_k + K_{kj}K^j_i
 - \Psi \left(K_{ik} + K_{ik}\right)\right.\vspace{2mm}\\ \hspace*{22mm}\left. - \Phi^{-2}\left(b_{i\vert k}
+ b_{k\vert i}\right)\right\} - \Phi^2b_ib_k,\vspace{2mm}\\
(b) \ \ \  R_{i00k} = \Theta_{ik\vert 0} - \Psi\Theta_{ik} + \Theta_{ij}\theta_k^j + \Phi^{4}\omega_{ij}\omega^j_k -
 \Phi^2b_ib_k \vspace{2mm}\\ \hspace*{22mm} - \frac{1}{2}\Phi^{2}\left(b_{i\vert k} +
b_{k\vert i}\right),\end{array}\eqno(4.9)$$
and

$$\begin{array}{lc}(a) \ \ \  K_{ik\vert 0} - K_{ki\vert 0} + K_{ij}K^j_k - K_{kj}K^j_i
 - \Psi \left(K_{ik} - K_{ki}\right)\vspace{2mm}\\ \hspace*{22mm} - \Phi^{2}\left(b_{i\vert k}
- b_{k\vert i}\right) = 0,\vspace{2mm}\\
(b) \ \ \  \omega_{ik\vert 0} = - \Psi\omega_{ik} + \Phi^{2}\left\{\omega_{kj}\Theta^j_i - \omega_{ij}\Theta^j_k
 + \frac{1}{2}\left(b_{i\vert k} - b_{k\vert i}\right)\right\}.\end{array}\eqno(4.10)$$
In particular, suppose that $(M, g)$ is a {\it vorticity-free spacetime}, that is, the vorticity tensor field
vanishes identically on $M$. Then, from (4.7c) and (4.8) we deduce that the curvature tensor field of the
Riemannian spatial connection satisfies the identities (4.7a), (4.7b) and the following:

$$\begin{array}{lc}(a) \ \ \ \bar{R}_{iljk} = \bar{R}_{jkil},\vspace{2mm}\\
(b) \ \ \ \bar{R}_{jkoi} =  \Theta_{ij\vert k} - \Theta_{ik\vert j} + \Theta_{ij}a_k -
\Theta_{ik}a_j. \end{array}\eqno(4.11)$$
Also, from (4.10b) we see that

$$b_{i\vert k} = b_{k\vert i},$$
and (4.6) becomes

$$\begin{array}{l}(a) \ \ \
R_{iljk} = \bar{R}_{iljk} + \Phi^{-2}\left(\Theta_{ij}\Theta_{lk} - \Theta_{ik}\Theta_{lj}\right),\vspace{2mm}\\
(b) \ \ \ R_{i0jk} = \Theta_{ik\vert j} - \Theta_{ij\vert k}+ \Theta_{ij}c_k - \Theta_{ik}c_j, \vspace{2mm}\\
(c) \ \ \ R_{il0k} = \bar{R}_{il0k} + \Theta_{lk}b_i - \Theta_{ik}b_l, \vspace{2mm}\\
(d) \ \ \  R_{i00k} = \Theta_{ik\vert 0} - \Psi\Theta_{ik} + \Theta_{ij}\Theta_{k}^j
- \Phi^{2}\left(b_{i\vert k} + b_ib_k\right).\end{array} \eqno(4.12)$$.

\section{Ricci tensor field and scalar curvature of a spacetime expressed in terms of spatial tensor fields}

Let $(M, g)$ be a $4D$ spacetime, and  $\{E_k\}$ be an orthonormal basis in $\Gamma(SM).$ Then, $\{\Phi^{-1}
\frac{\partial}{\partial x^0}, E_k\}$ is an orthonormal frame field on $M$. According to  \cite{o}, p.87, the
Ricci tensor of $(M, g)$ is given by

$$Ric(X, Y) = \sum_{k =1}^3R(E_k, X, E_k, Y) - \Phi^{-2}R(\frac{\partial }{\partial x^0}, X,
\frac{\partial }{\partial x^0}, Y), \eqno(5.1)$$
for all $X, Y \in\Gamma(TM).$ Now, we express $E_k$ as follows

$$E_k = E_k^i\frac{\delta }{\delta x^i},\eqno(5.2)$$
and obtain

$$\bar{g}^{ij} =  \sum_{k =1}^3E_k^iE_k^j.\eqno(5.3)$$
Then we consider the following local components of{\it Ric} with respect to the threading frame
$\{\partial/\partial x^0, \delta/\delta x^i \}:$

 $$\begin{array}{lc}(a) \ \ \  R_{ij} =
 Ric\left(\frac{\delta }{\delta x^j}, \frac{\delta }{\delta x^i}\right), \ \ \ \ (b) \ \ \ R_{i0} =
 Ric\left(\frac{\partial }{\partial x^0}, \frac{\delta }{\delta  x^i}\right),
\vspace{4mm}\\ \ \ (c) \ \ \  R_{00} =  Ric\left(\frac{\partial }{\partial x^0},
\frac{\partial }{\partial x^0}\right),\end{array}\eqno(5.4) $$
and by  using (5.2), (5.3), (5.4) and $(A_9)$ into (5.1), we obtain

$$\begin{array}{c}(a) \ \
 R_{ij} =  \bar{g}^{lk}R_{iljk} +  \Phi^{-2}R_{i00j}, \ \ \ (b) \ \ R_{i 0} = \bar{g}^{lk}R_{il0k}
= -\bar{g}^{lk}R_{k0il},
\vspace{2mm}\\ \ \ (c) \ \ \ R_{00} =  -\bar{g}^{ik}R_{i00k}.\end{array}\eqno(5.5) $$
Now, by using (4.5), (4.6) and (4.9) in (5.5), we deduce that

$$\begin{array}{lc} (a) \ \ \ R_{ij} = \bar{R}_{i\ \ jk}^{\;k} + \Phi^{-2}\left\{\left(\Theta -
 \Psi\right)K_{ij} + K_{ij\vert 0}\right\} - b_{i\vert j} - b_ib_j,  \vspace{2mm}\\
 (b) \ \ \ R_{i0} = \bar{R}_{i\ 0k}^{\;k} + \Theta b_i - K_{ik}b^k\vspace{2mm}\\ \hspace*{22mm}
 = K^k_{i\vert k} - \Theta_{\vert i} + \Theta c_i - c_kK_I^k + 2\Phi^2\omega_{ik}b^k,\vspace{2mm}\\
(c) \ \ \ R_{00} = \Psi\Theta - \Theta_{\vert 0} - K_h^kK_k^h  +
\Phi^{2}\left\{b^2 + b^k_{\ \vert k}\right\}, \end{array} \eqno(5.6)$$
and

$$\begin{array}{l} (a) \ \ \ R_{ij} = \bar{R}_{i\ \ jk}^{\; k} + \Phi^{-2}\left\{\left(\Theta - \Psi\right)
 \Theta_{ij} + \Theta_{ij\vert 0}\right\} + \omega_{ij\vert 0} \vspace{2mm}\\
\hspace*{22mm} + \left(\Theta + \Psi\right)\omega_{ij} - b_{i\vert j} - b_ib_j,  \vspace{2mm}\\
(b) \ \ \ R_{i0} = \bar{R}_{i\ \ 0k}^{\; k} + \Theta b_{i} - \Theta_{ik}b^k - \Phi^{2}\omega_{ik}b^k,\vspace{2mm}\\
\hspace*{22mm}= \Theta_{i\vert k}^k - \Theta_{\vert i} + \Theta c_i - \Theta_{ik}c^k
\vspace{2mm}\\ \hspace*{22mm} - \Phi^2\left(\omega_{i\vert k}^k  + c_k\omega_i^k  - \omega_{ik}b^k\right),
 \vspace{2mm}\\
(c) \ \ \ R_{00} = \Psi\Theta - \Theta_{\vert_0} - \sigma^2 - \frac{1}{3}\Theta^2 + \Phi^4\omega^2 +
\Phi^2\left(b^2 + b^k_{\vert k}\right), \end{array} \eqno(5.7)$$
where we put

$$\sigma^2 = \sigma_{hk}\sigma^{hk}, \ \ \ \omega^2 = \omega_{hk}\omega^{hk}, \ \ \ b^2 = b_{k}b^{k}.$$
Next, we take the symmetric and skew-symmetric parts in (5.6a) and (5.7a), and obtain

$$\begin{array}{lc}(a) \ \  \ R_{ij} = \bar{R}_{ij} + \Phi^{-2}\left\{\left(\Theta - \Psi\right)\Theta_{ij} +
 \Theta_{ij\vert 0}\right\}\vspace{2mm}\\ \hspace*{22mm} - \frac{1}{2}\left(b_{i\vert j} +
b_{j\vert i}\right) - b_ib_j, \vspace{3mm}\\
(b) \ \ \ \frac{1}{2}\left(\bar{R}_{i\ jk}^{\;k} -\bar{R}_{j\ ik}^{\;k}\right) = \frac{1}{2}\left(b_{i\vert j} -
b_{j\vert i}\right) - \omega_{ij\vert 0} - \left( \Theta + \Psi\right)\omega_{ij}, \end{array}\eqno(5.8)$$
where we put

 $$\bar{R}_{ij} = \frac{1}{2}\left(\bar{R}_{i\ jk}^{\;\;k} + \bar{R}_{j\ ik}^{\;\;k}\right). \eqno(5.9)$$
The spatial with local components $\bar{R}_{ij}$ is called the {\it spatial Ricci tensor} of the  spacetime
$(M, g)$.\par
 The above formulas for Ricci tensor enable us to obtain a new formula for the scalar curvature ${\bf R}$ of $(M, g)$
 in terms of spatial tensor fields. We start with ${\bf R}$ given by

 $${\bf R} = \sum_{k=1}^3Ric(E_k, E_k) - \Phi^{-2}Ric\left(\frac{\partial }{\partial x^0},
\frac{\partial }{\partial x^0}\right) = \bar{g}^{ij}R_{ij} - \Phi^{-2}R_{00}.\eqno(5.10)$$
Then by direct calculations, using (5.8a) and (5.7c) into (5.10), we obtain

$${\bf R} = \bar{{\bf R}} + \Phi^{-2}\left\{\frac{4}{3}\Theta^2 - 2\Psi\Theta + 2\Theta_{\vert 0}\right\}
- \Phi^2\omega^2 - 2b^2 - 2b^k_{\ \vert k},\eqno(5.11)$$
where we put

$$\bar{{\bf R}} = \bar{g}^{ij}\bar{R}_{ij}.\eqno(5.12) $$
We call $\bar{{\bf R}}$ the {\it spatial scalar curvature} of the spacetime $(M, g)$. Note that both $\bar{R}_{ij}$
and $\bar{{\bf R}}$ are related to the geometry of the spatial distribution $SM$ which is not necessarily
supposed to be integrable.

\section{The splitting of both the Einstein gravitational tensor field and the energy-momentum tensor field}

We start with the Einstein gravitational tensor field $G$ of $(M, g)$ given by

$$G(X, Y) = Ric(X, Y) - \frac{\bf R}{2}g(X, Y), \ \ \ \forall \ X, Y\ \in \Gamma(TM).\eqno(6.1)$$
Then, with respect to the threading frame field $\{\partial/\partial x^0, \delta/\delta x^i\}$ we have

$$\begin{array}{lc}(a) \ \ \ G_{ij} = G(\frac{\delta }{\delta x^j}, \frac{\delta }{\delta x^i}) = R_{ij} -
\frac{\bf R}{2}g_{ij},\vspace{2mm}\\  (b) \ \ \ G_{i0} = G(\frac{\partial }{\partial x^0},
\frac{\delta }{\delta x^i}) = R_{i0},\vspace{2mm}\\ (c) \ \ \ G_{00} = G(\frac{\partial }{\partial x^0},
 \frac{\partial }{\partial x^0}) = R_{00} + \Phi^2\frac{\bf R}{2}.\end{array}\eqno(6.2)$$
By using (5.8a), (5.7b), (5.7c) and (5.11) into (6.2), we obtain

$$\begin{array}{lc}(a) \ \ \ G_{ij} = \bar{G}_{ij} + \left(b^2 + b_{\;\;\vert k}^k + \frac{1}{2}\Phi^2\omega^2
\right)\bar{g}_{ij} - \frac{1}{2}\left(b_{i\vert j} +  b_{j\vert i}\right) - b_ib_j \vspace{2mm}\\ \hspace*{8mm}
+ \Phi^{-2}\left\{\left( \Theta - \Psi\right)\Theta_{ij} + \Theta_{ij\vert 0} - \left(\Theta_{\vert 0} - \Psi\Theta
+ \frac{2}{3}\Theta^2 + \frac{1}{2}\sigma^2\right)\bar{g}_{ij}\right\}, \vspace{2mm}\\
(b) \ \ \ G_{i0} = \bar{R}_{i\ 0k}^{\;\;k} + \Theta b_i - \Theta_{ik}b^k - \Phi^2\omega_{ik} \vspace{2mm}\\
\hspace*{12mm} = \Theta^k_{i\vert k} - \Theta_{\vert i} + \Theta c_i - \Theta_{ik}c^k - \Phi^2\left(
\omega^k_{i\vert k} + c_k\omega^k_{i} - 2\omega_{ik}b^k\right),\vspace{2mm}\\
(c)\ \ \ G_{00} = \frac{1}{2}\left(\Phi^2\bar{\bf R} + \frac{2}{3}\Theta^2 - \sigma^2 +\Phi^4\omega^2\right),
\end{array}\eqno(6.3)$$
where we put

$$\bar{G}_{ij} = \bar{R}_{ij} - \frac{\bar{\bf R}}{2}\bar{g}_{ij}.\eqno(6.4)$$
The spatial tensor field $\bar{G}$ with local components $\bar{G}_{ij}$ is called the {\it $3D$ Einstein
gravitational tensor field} of the spacetime $(M, g)$.\par
Now, in order to give a coordinate-free formula for the energy-momentum tensor field, we consider a spatial
1-form $\omega$ and a spatial tensor field $S$ of type $(0, 2)$. Then, we define a 1-form and a
tensor field of type $(0, 2)$ on $M$, denoted by the same symbols and given by

$$(a) \ \ \ \ \omega(X) = \omega({\cal{S}}X), \ \ \ \ (b) \ \ \ \ S(X, Y) = S({\cal{S}}X, {\cal{S}}Y),\eqno(6.5)$$
for all $X, Y \ \in \ \Gamma(TM)$. As an example, the Riemannian metric $\bar{g}$ on $SM$ defines a symmetric
tensor field $\bar{g}$ on $M$ given by

$$ \bar{g}(X, Y) = \bar{g}({\cal{S}}X, {\cal{S}}Y),\ \ \ \forall \ X, Y \ \in \ \Gamma(TM).\eqno(6.6)$$
Note that $\bar{g}$ from (6.6) coincides with the tensor field $h$ given by its local components in formula (4.10)
of \cite{ej}. Also, we need the 1-form $u$ induced by the unit vector field $U = \Phi^{-1}\partial/\partial x^0$
by the formula

$$u(X) = \Phi^{-1}g\left(X, \frac{\partial }{\partial x^0}\right), \forall \ X \in \ \Gamma(TM).\eqno(6.7)$$
Based on these geometric objects, we claim that the energy-momentum tensor field $T$ measured by an observer
moving with the unit 4-velocity $U$ has the following coordinate-free expression:

$$\begin{array}{lc}T(X, Y) = \rho u(X)u(Y) + q(X)u(Y) + q(Y)u(X) + p\bar{g}(X, Y)\vspace{2mm} \\
\hfill  + \ \pi(X, Y), \forall \ X, Y \ \in \ \Gamma(TM).\end{array}\eqno(6.8)$$
Here, $\rho$ and $p$ are the relativistic energy density and the relativistic pressure respectively, while $q$
is a 1-form on $M$ defined by a spatial 1-form as in (6.5a) and $\pi$ is a symmetric and trace-free tensor
field on $M$ defined by a spatial tensor field as in (6.5b). Now, take $X = Y = U$ in (6.8) and obtain

$$\rho = T(U, U).\eqno(6.9)$$
Then we put

$$T_{00} = T(\frac{\partial }{\partial x^0}, \frac{\partial }{\partial x^0}), \eqno(6.10)$$
and (6.9) becomes

$$T_{00} = \Phi^2\rho.\eqno(6.11)$$
Similarly, take $X = U$ and $Y = {\cal{S}}Y$ in (6.8) and deduce that

$$q({\cal{S}}Y) = -T(U, {\cal{S}}Y).\eqno(6.12)$$
Now, we put

$$(a) \ \ \ \ T_{i0} = T(\frac{\partial }{\partial x^0}, \frac{\delta }{\delta x^i}), \ \ \ \ (b) \ \ \ \
q_i = q(\frac{\delta }{\delta x^i}),\eqno(6.13)$$
and taking $Y = \delta /\delta x^i$ in (6.12), we infer that

$$T_{i0} = -\Phi q_i.\eqno(6.14)$$
Finally, denote

$$(a) \ \ \ \ T_{ij} = T(\frac{\delta }{\delta x^j}, \frac{\delta }{\delta x^i}), \ \ \ \ (b) \ \ \ \
\pi_{ij} = \pi(\frac{\delta }{\delta x^j}, \frac{\delta }{\delta x^i}),\eqno(6.15)$$
and taking $X = \delta /\delta x^j$, $Y = \delta /\delta x^i$ in (6.8), obtain

$$T_{ij} = p\bar{g}_{ij} + \pi_{ij}.\eqno(6.16)$$
Contracting (6.16) by $\bar{g}^{ij}$ and taking into account that $\pi$ is defined by a trace-free spatial
tensor field, we infer that

$$(a) \ \ \ \ p = \frac{1}{3}\bar{g}^{ij}T_{ij}, \ \ \ \ (b) \ \ \ \ \pi_{ij} = T_{ij} -
\frac{1}{3}\left(T_{hk}\bar{g}^{hk}\right)\bar{g}_{ij}.\eqno(6.17)$$
Now, taking into account(6.9) - (6.17) it is clear that (6.8) represents the coordinate-free version of (5.9)
from \cite{ej}, p.91. Thus $q$ and $\pi$ from (6.8) are the {\it relativistic momentum density} and the {\it
relativistic anisotropic (trace-free) stress tensor field}, respectively.

\section{A new splitting of Einstein field equations with respect to a general $(1+3)$ threading of spacetime}

We start this section with the coordinate-free form of Einstein field equations $(EFE)$ expressed as follows:

$$G(X, Y) + \Lambda g(X, Y) = 8\pi {\bf G}T(X, Y), \ \ \ \forall \ X, Y \in \Gamma(TM), \eqno(7.1)$$
where $\Lambda$ is the cosmological constant and ${\bf G}$ is the Newton constant. Now, take
$X = \delta /\delta x^j$, $Y = \delta /\delta x^i$ in (7.1) and by using (6.3a), (2.9) and (6.16), we obtain

$$\begin{array}{lc}\bar{G}_{ij} + \left(b^2 + b_{\;\;\vert k}^k + \frac{1}{2}\Phi^2\omega^2 + \Lambda
\right)\bar{g}_{ij} - \frac{1}{2}\left(b_{i\vert j} +  b_{j\vert i}\right) - b_ib_j \vspace{2mm}\\
+ \Phi^{-2}\left\{\left( \Theta - \Psi\right)\Theta_{ij} + \Theta_{ij\vert 0} - \left(\Theta_{\vert 0} -
\Psi\Theta + \frac{2}{3}\Theta^2 + \frac{1}{2}\sigma^2\right)\bar{g}_{ij}\right\}\vspace{2mm} \\ = 8\pi {\bf G}
\left(p\bar{g}_{ij} + \pi_{ij}\right).\end{array}\eqno(7.2)$$
The equations from (7.2) will be called the {\it spatial Einstein field equations} (SEFE). Next, we take
$X = \partial /\partial x^0$, $Y = \delta /\delta x^i$ in (7.1) and by using (6.3a) and (6.14), we deduce that

$$ \begin{array}{lc}(a) \ \ \ \
\bar{R}_{i\ 0k}^{\;\;k} + \Theta b_i - \Theta_{ik}b^k - \Phi^2\omega_{ik} = -8\pi{\bf G}\Phi q_i, \vspace{2mm}\\
(b) \ \ \ \ \Theta^k_{i\vert k} - \Theta_{\vert i} + \Theta c_i - \Theta_{ik}c^k - \Phi^2\left(
\omega^k_{i\vert k} + c_k\omega^k_{i} - 2\omega_{ik}b^k\right)\vspace{2mm}\\ \hspace*{22mm}
= -8\pi{\bf G}\Phi q_i.\end{array}\eqno(7.3)$$
The equations from either (7.3a) or (7.3b) are called {\it mixed Einstein field equations} $(MEFE)$. Finally,
take $X = Y = \partial /\partial x^0$ in (7.1), and by using (6.3c), (2.2a) and (6.11), we infer that

$$\Phi^2\bar{\bf R} + \frac{2}{3}\Theta^2 - \sigma^2 +\Phi^4\omega^2 = 2\Phi^2\left(\Lambda +
8\pi{\bf G}\rho\right).\eqno(7.4)$$
We call (7.4) the {\it temporal Einstein field equation} $(TEFE)$. Thus the Einstein field equations (7.1)
are splitting in three groups of equations given by (7.2), (7.3) and (7.4). {\it It is worth mentioning that
these equations are expressed in terms of spatial tensor fields and their covariant derivatives induced by
the Riemannian spatial connection}.\par
Next, by contracting (7.2) by $\bar{g}^{ij}$ and using (6.4) we deduce that the spatial scalar curvature is given by

$$\begin{array}{lc}\frac{1}{2}\bar{\bf R} = 2\left(b^2 + b^k_{\ \vert k}\right) + \frac{3}{2}\Phi^2\omega^2 +
3\Lambda \vspace{2mm}\\ \hspace*{8mm}  - \Phi^{-2}\left\{\Theta^2 + \Theta_{\vert 0} + \frac{3}{2}\sigma^2
- 2\Psi\Theta\right\}- 24\pi{\bf G}p. \end{array}\eqno(7.5)$$
Comparing (7.5) with (7.4) we obtain the {\it Raychaudhuri-Ehlers equation} induced by the general $(1+3)$
threading of the spacetime:

$$\Theta_{\vert 0} + \frac{1}{3}\Theta^2 + \sigma^2 - \Psi\Theta - \Phi^2\left\{b^2 + b^k_{\ \vert k} +
\Phi^2\omega^2 + \Lambda - 4\pi G(\rho + 3p)\right\} = 0.\eqno(7.6)$$
Note that (7.6) is the generalization of the equation (6.4) from \cite{ej}, which was obtained for the
particular case $\Phi = 1$.

\section{A new splitting of conservation laws with respect to a general $(1+3)$ threading of spacetime}

As is well known, the energy-momentum conservation equations are given by the vanishing of the divergence of $T$.
In order to obtain their explicit form, we consider an orthonormal frame field $\{E_k, U = \Phi^{-1}\partial/
\partial x^0\},$ and according to \cite{o}, p.86 we have

$$\left(div T\right)(X) = \sum_{k=1}^3(\nabla_{E_k}T)(E_k, X) - (\nabla_{U}T)(U, X) = 0, \eqno(8.1)$$
for all $X \in \Gamma(TM)$. Then, take $X = \partial/\partial x^0,$ in (8.1), and by using $(B_1)$ and $(B_2)$
 obtain

$$\Phi^2\left\{T^j_{0\vert j} + \left(2b^j - c^j\right)T_{j0} - \Theta^{jk}T_{jk}\right\} - \frac{\partial T_{00}}
{\partial x^0} + \left(2\Psi - \Theta\right)T_{00} = 0.\eqno(8.2)$$
Similarly, take $X = \delta/\delta x^i$ in (8.1) and by using $(B_3)$ and $(B_4)$ we deduce that

$$\Phi^2\left\{T^j_{i\vert j} - \omega_{ij}T^j_0 + T_{ij}b^j\right\} + \left(\Psi - \Theta\right)T_{i0} -
T_{i0\vert 0} + T_{00}b_i - \Theta_{ij}T^j_0 = 0.\eqno(8.3)$$
Now, from (6.14) and (6.16), we infer that

$$(a) \ \ \ T^j_0 = -\Phi q^j, \ \ \ \ (b) \ \ \ \ T_i^j = p\delta_i^j + \pi_i^j.\eqno(8.4)$$
Then, by using (2.19), (6.11), (6.14), (6.16) and (8.4) in (8.2) and (8.3), and taking into account that
both $\sigma_{ij}$ and $\pi_{ij}$ are trace-free spatial tensor fields, we obtain the {\it energy conservation
 equation}:

$$\frac{\partial \rho}{\partial x^0} + \left(\rho + p\right)\Theta + \sigma^{jk}\pi_{jk} + \Phi\left(q^j_{\ \vert j}
+ 2b_jq^j\right) = 0,\eqno(8.5)$$
and the {\it momentum conservation equation}:

$$q_{i\vert 0} + \frac{4}{3}\Theta q_i + \left(\sigma_{ij} + \Phi^2\omega_{ij}\right)q^j + \Phi
\left\{\frac{\delta p}{\delta x^i} + \pi^j_{i\vert j} + \left(p + \rho\right)b_i + \pi_{ij}b^j\right\} = 0.\eqno(8.6)$$
In order to compare with what is known in the literature with respect to the $(1+3)$ threading of spacetime,
we note that

$$\nabla_UU = b^k\frac{\delta }{\delta x^k}.\eqno(8.7)$$
This is obtained by direct calculations using (2.18d). Also, we should remark that (8.7) states that though the
velocity $U$ is a timelike vector field, the acceleration $\nabla_UU$ is a spatial vector field.\par
{\bf Remark 8.1} The above conservation laws are obtained in the most general setting. Indeed, if in particular
$\Phi = 1$, (8.5) and (8.6) become (5.11) and (5.12) from \cite{ej}, p.92, respectively. If moreover, we have
a perfect fluid on $M$, that is, $q_i = 0$ and $\pi_{ij} = 0$, then (8.5) and (8.6) become (1) and (2)
from Proposition 5 in \cite{o}, p.339 .

\section{Splitting of Einstein field equations in an almost $FLRW$ universe}

Let $M$ be an $FLRW$ universe, whose line element is given by

$$d\tilde{s}^2 = \tilde{g}_{ab}dx^adx^b = a^2(\tau)\left\{-d\tau^2 + \delta_{ij}dx^idx^j\right\}, \eqno(9.1)$$
where $\tau$ is the conformal time on $M$ and the three-dimensional space given by $\tau$ = const., is an Euclidean
space. The Lorentz metric given by (9.1) is called the {\it background metric}. For a more realistic model of
the universe there have been studied perturbation of this metric (cf.\cite{ej}, \cite{mw}, \cite{mf}). The line
element of the {\it full metric} $g$ on $M$ is expressed as follows:

$$ds^2 = g_{ab}dx^adx^b = \tilde{g}_{ab}dx^adx^b + \delta g_{ab}dx^adx^b, \eqno(9.2)$$
where $\delta g_{ab}$ determine the perturbation. In the present paper we consider the conformal-Newtonian gauge
case, for which the full metric is given by

$$ds^2 = a^2(\tau)\left\{-(1 + 2A)d\tau^2 + (1 - 2B)\delta_{ij}dx^idx^j\right\}, \eqno(9.3)$$
where $A$ and $B$ are the well-known {\it Bardeen invariants}. The spacetime $(M, g)$ with $g$ given by (9.3)
 is called an {\it almost $FLRW$ universe}.\par
Now, by using (2.2b) for (9.3) we deduce that $\xi_i = 0, \ i \in\{1,2,3\},$ which imply

$$\begin{array}{lc}(a) \ \ \ \frac{\delta }{\delta x^i} = \frac{\partial }{\partial x^i}, \ \ \ \ (b) \ \ \ \
\omega_{ij} = 0, \ \ \ \ (c) \ \ \ \ a_i = 0,\vspace{2mm} \\ (d) \ \ \ \ \bar{g}_{ij} = a^2(1 - 2B)\delta_{ij} \ \ \ \
(e) \ \ \ \ \bar{g}^{ij} = \frac{1}{a^2(1 - 2B)}\delta^{ij}.\end{array}\eqno(9.4)$$
Also, according to the notation in (2.2a), we have
$$\Phi^2 = a^2(1 + 2A).\eqno(9.5)$$

From (9.4b) we see that the spatial distribution $SM$ of $(M, g)$ is integrable, but its leaves are not
anymore Euclidean spaces. Moreover, by using (2.11a), (2.11b), (2.6d), (9.4d) and (9.5), we obtain

$$\begin{array}{lc}(a) \ \ \ \Theta_{ij} = a^2\left\{(1-2B){\cal{H}} - B^{\prime}\right\}\delta_{ij}, \ \ \ (b)
 \ \ \ \Theta = 3\left\{{\cal{H}} - \frac{B^{\prime}}{1 - 2B}\right\}, \vspace{2mm} \\ (c) \ \ \ \Psi = {\cal{H}} +
 \frac{A^{\prime}}{1 + 2A}, \end{array}\eqno(9.6)$$
where $"\vert "$ denotes derivative with respect to $\tau$, and ${\cal{H}} = a^{\prime}/a$ is the {\it Hubble
 parameter} of the background metric. By calculations, using (9.6a), (9.6b) and (9.4d) into (2.11c),
we infer that

$$ (a) \ \ \ \sigma_{ij} = 0, \ \ \ (b) \ \ \ \Theta_{ij} = \frac{1}{3}\Theta\bar{g}_{ij}, \ \ \forall \ i, j
\in \{1, 2, 3\}. \eqno(9.7)$$
Taking into account that $\Phi^{-2}\Theta_{ij}$ are the local components of the second fundamental form of the
leaves of $SM$ (see (2.18a)), and using (9.7), we can state the following result on the kinematic quantities
and geometry of $(M, g)$.\vspace{2mm} \par

{\bf Theorem 9.1} {\it Let $(M, g)$ be an almost $FLRW$ universe. Then we have the following assertions}:\par
(i) $(M, g)$ {\it is both vorticity-free and shear-free spacetime}.\par
(ii) {\it The leaves of the spatial distribution are totally umbilical hypersurfaces of $(M, g)$ with mean
curvature vector}

$$H = \frac{1}{3}\Phi^{-2}\Theta\frac{\partial }{\partial\tau}.\eqno(9.8)$$

Next, we assume that the energy-momentum tensor for the almost $FLRW$ universe $(M, g)$ takes the perfect fluid
form, that is,

$$ q_i = 0, \ \ \ \pi_{ij} = 0, \ \ \forall \ i,j \in \{1,2,3\}.\eqno(9.9)$$
Then, the Bardeen invariants coincide (cf. \cite{ej}, p.259), that is, from now on in our calculations we put
$A$ instead of $B$. First, by using (2.19), (9.4c), (2.6b), (9.5) and (9.4e), we deduce that

$$(a) \ \ \ b_i = c_i = \frac{A_i}{1+2A}, \ \ \ (b) \ \ \ b^2 = \frac{1}{a^2(1-4A^2)(1+2A)}
\sum_{k=1}^3(A_k)^2, \eqno(9.10)$$
where we put $A_i = \partial A/\partial x^i$. The local coefficients of the spatial Riemannian connection
$\bar{\nabla}$ are given by ( see (2.15), $(C_1)$, (2.12), (9.7))

$$\begin{array}{lc}(a) \ \ \ \bar{\Gamma}^{\;\;k}_{i\ j} = \frac{1}{1-2A}\left\{\delta_{ij}A^k - \delta_i^kA_j -
\delta_j^kA_i \right\}, \vspace{2mm}\\ (b) \ \ \ K_i^j = \Theta_i^j = \frac{1}{3}\Theta\delta_i^j,
\end{array}\eqno(9.11)$$
where we have $A^k = A_k$. Then, by using (9.10a) and (9.11a), we obtain the spatial covariant derivative
and the divergence of the acceleration, given by

$$\begin{array}{lc}(a) \ \ \ b_{i\vert j} = \frac{1}{1 + 2A}\left\{A_{ij} + \frac{8A}{1 - 4A^2}A_iA_j -
\frac{1}{1 - 2A}\delta_{ij}\displaystyle{\sum_{k=1}^3}(A_k)^2\right\},\vspace{2mm}\\
(ba) \ \ \ b^k_{\ \vert k} = \frac{1}{a^2(1 - 4A^2)}\displaystyle{\sum_{k=1}^3}\left\{A_{kk} + \frac{2A - 3}{1 - 4A^2}
(A_k)^2\right\}, \end{array}\eqno(9.12)$$
where we put

$$A_{ij} = \frac{\partial^2A}{\partial x^i\partial x^j}.$$
Now, by direct calculation using (9.6b), we infer that

$$\begin{array}{lc}(a) \ \ \ \Theta_{\vert 0} = \frac{\partial \Theta}{\partial\tau} = 3\left({\cal{H}}^{\prime} -
\frac{2(A^{\prime})^2 + (1-2A)A^{\prime\prime}}{(1 - 2A)^2}\right),\vspace{2mm}\\
(b) \ \ \ \Theta_{\vert i} = \frac{\partial \Theta}{\partial x^i} = -\frac{3}{(1 - 2A)^2}\left\{2A^{\prime}A_i +
(1 - 2A)A^{\prime}_i\right\}.\end{array}\eqno(9.13)$$
Also, by using (9.7b), (9.11b) and (9.13), and taking into account that $\bar{\nabla}$ is a metric connection, we obtain

$$\begin{array}{lc}(a) \ \ \ \Theta_{ij\vert 0} = \frac{1}{3}\Theta_{\vert 0}\bar{g}_{ij} \vspace{2mm}\\
 \hspace*{18mm} =
\frac{a^2}{1 - 2A}\left\{(1 - 2A)^2{\cal{H}}^{\prime} - (1-2A)A^{\prime\prime}
 -2(A^{\prime})^2\right\}\delta_{ij},\vspace{2mm}\\
(b) \ \ \ \Theta^k_{i\vert k} = \frac{1}{3}\Theta_{\vert i} = -\frac{1}{(1 - 2A)^2}\left\{2A^{\prime}A_i +
(1 - 2A)A^{\prime}_i\right\}.\end{array}\eqno(9.14)$$
Finally, by using $(C_6)$ and $(C_7)$ in (6.4), we deduce that the $3D$ Einstein gravitational tensor field is given by

$$\bar{G}_{ij} = \frac{1}{(1 - 2A)^2}\left\{3A_iA_j + (1 - 2A)A_{ij} - \delta_{ij}\sum_{k=1}^3\left(2(A_k)^2 +
(1 - 2A)A_{kk}\right)\right\}.\eqno(9.15)$$
Now, we are in a position to present the splitting of $(EFE)$ for the almost universe $(M, g)$. First, we consider
that the background energy-momentum tensor $\tilde{T}$ take the perfect fluid form, that is, we have
$$(a) \ \ \ \tilde{T}_{ij} = \tilde{p}\tilde{g}_{ij} = \tilde{p}a^2\delta_{ij}, \ \ \ (b) \ \ \ \tilde{T}_{io} = 0,
\ \ \ (c) \ \ \ \tilde{T}_{00} = a^2\tilde{\rho}, \eqno(9.16)$$
where $\tilde{\rho}$ and $\tilde{p}$ are the relativistic density and pressure, respectively. Then, the energy
momentum tensor $T$ of the perturbed universe $(M, g)$ should have the perfect fluid from too. Hence, by (6.16),
(6.14) and (6.11), we have

$$\begin{array}{lc}(a) \ \ \ T_{ij} = p\bar{g}_{ij} = pa^2(1 - 2A)\delta_{ij}, \ \ \
(b) \ \ \ T_{io} = 0,\vspace{2mm} \\ \hspace*{22mm} (c) \ \ \ T_{00} = \Phi^2\rho = a^2(1 + 2A)\rho,
\end{array}\eqno(9.17)$$
where we put

$$p = \tilde{p} + \delta p\ \ \ \ \mbox{and} \ \ \ \ \rho = \tilde{\rho} + \delta\rho.$$
After some long calculations using (9.15), (9.10), (9.12), (9.13a), (9.6), (9.7), (9.4d), (9.5) and (9.17a) into
(7.2) we deduce that the $(SEFE)$ for the almost $FLRW$ universe $(M, g)$ are given by

$$\begin{array}{lc}\frac{2}{1 - 4A^2}\left\{2AA_{ij} + \frac{1 + 4A + 12A^2}{1 - 4A^2}A_iA_j \right\}\vspace{2mm}\\
+ \delta_{ij}\left\{\frac{1}{1 + 2A}\left(2A^{\prime\prime} + \frac{6 + 4A}{1 + 2A}{\cal{H}}A^{\prime} +
\frac{6A - 1}{1 - 4A^2}(A^{\prime})^2 - (1 -2A)({\cal{H}}^2 + 2{\cal{H}}^{\prime})\right)\right.\vspace{2mm}\\
\left.- \frac{1}{1 - 4A^2}
\left(4A\displaystyle{\sum_{k=1}^3}A_{kk} + \frac{3 + 4A + 12A^2}{1 - 4A^2}\displaystyle{\sum_{k=1}^3}
 (A_{k})^2\right)\right\}\vspace{2mm}\\
 = \left(8\pi {\bf G}p - \Lambda\right)a^2(1 - 2A)\delta_{ij}.\end{array}\eqno(9.18)$$
Next, since $(M, g)$ is vorticity-free and $q_i = 0$, (7.3b) becomes

$$\Theta_{\vert i} - \Theta c_i = 0,\eqno(9.19)$$
via (9.7b) and (9.14b). Then, by using (9.6b), (9.10a) and (9.13b) into (9.19), we obtain the following $(MEFE)$ for
$(M, g)$:

$$(1 + 6A)A^{\prime}A_i + (1 - 4A^2)A_i^{\prime} + (1 - 2A)^2{\cal{H}}A_i = 0.\eqno(9.20)$$
Also, by using $(C_7)$, (9.6b) and (9.5) into (7.4), and taking into account that $\sigma^2 = \omega^2 = 0,$ we
deduce that the $(TEFE)$ for $(M, g)$ is given by

$$\begin{array}{lc}\frac{1 + 2A}{(1 - 2A)^3}\displaystyle{\sum_{k=1}^3}\left\{3(A_k)^2 + 2(1 - 2A)A_{kk}\right\}
 + 3\left({\cal{H}} - \frac{A^{\prime}}{1 - 2A}\right)^2\vspace{2mm}\\ \hspace*{22mm} = a^2(1 + 2A)\left(\Lambda
+ 8\pi{\bf G}\rho\right).\end{array}\eqno(9.21)$$
Finally, the Raychaudhuri - Ehlers equation (7.6) becomes

$$\begin{array}{lc}3\left\{{\cal{H}}^{\prime} - \frac{1}{1 - 2A}A^{\prime\prime} - \frac{1}
{1 - 4A^2}{\cal{H}}A^{\prime} - \frac{4A}{(1 - 2A)(1 - 4A^2)}(A^{\prime})^2\right\}\vspace{2mm}\\ -
\frac{1}{1 - 2A}\left\{\displaystyle{\sum_{k=1}^3}A_{kk} -
\frac{2}{1 - 4A^2}\displaystyle{\sum_{k=1}^3}(A_k)^2\right\}\vspace{2mm}\\ =
a^2(1 + 2A)\left(\Lambda - 4\pi{\bf G}(\rho + 3p)\right).\end{array}\eqno(9.22)$$
Summing up the above results, we state the following theorem.\vspace{2mm}\newline
{\bf Theorem 9.2} {\it Let $(M, g)$ be an almost $(FLRW)$ universe with the energy-momentum tensor of a perfect
fluid. Then we have the following assertions}:\par
(i) {\it The $(SEFE)$, $(MEFE)$ and $(TEFE)$ of $(M, g)$ are given by} (9.18), (9.20) {\it and} (9.21),
{\it respectively}.\par
(ii) {\it The Raychaudhuri-Ehlers equation in $(M, g)$  is given by} (9.22).\vspace{2mm}\par
In particular, suppose that $A = 0$ and $\Lambda = 0$. Then, (9.21) and (9.22) become

$${\cal{H}}^2 = \frac{8\pi{\bf G}}{3}a^2\rho, \eqno(9.23)$$
and

$${\cal{H}}^{\prime} = - \frac{4\pi{\bf G}}{3}a^2(\rho + 3p), \eqno(9.24)$$
respectively. Note that in this case, (9.20) is trivial, and (9.18) is a consequence of (9.23) and (9.24). As
(9.23) and (9.24) are the well known Friedmann equations for an $(FLRW)$ universe (cf. (8.4), and (8.5) in
\cite{mw}), we are entitled to call (9.21) and (9.22) the {\it perturbed Friedmann equations} in an almost
$(FLRW)$ universe.

\section{Conclusions}

The idea to develop a $(1+3)$ threading of a spacetime with respect to a non-normalized vector field came up as
a need for the study of the spacetimes whose metrics have the general form (2.2) with $\Phi^2 \neq 1$. We only
mention here all the metrics from both the theory of cosmological perturbations and the theory of black holes.\par
The main difference between our approach and the methods developed in earlier papers consists in the fact that we
deal with the spatial tensor fields as intrinsic objects from the geometry of the spatial distribution. In
earlier literature, the spatial tensor fields have been considered as projections on the spatial distribution
of tensor fields defined on $M$. This was the main obstacle in defining a correct spatial covariant
differentiation for the general case when the spatial distribution is not integrable (cf. (4.19) of
\cite{ej}).\par
The approach we develop in the paper is based on spatial tensor fields and on the Riemannian spatial connection
which behave as $3D$ geometric objects with local components defined on a $4D$ spacetime. It is noteworthy that
the three groups of $(EFE)$ presented in the paper are expressed in terms of spatial tensor fields and their
covariant derivatives induced by the Riemannian spatial connection. This enable
 us to write down a splitting of $(EFE)$ for an almost $FLRW$ universe, which might have an important role in
the difficult task of finding models for such a universe. Moreover, the approach can be extended to the study of
threading of higher-dimensional universes. This can be seen in a paper of the first author on the threading of a
$5D$ universe (cf. \cite{b}).

\section*{Appendix A}

We shall present here details about the formulas stated in Section 4. For the curvature tensor $R$ of
the Levi-Civita connection $\nabla$ on $(M, g)$ we use the formula

$$R(X, Y)Z = \nabla_X\nabla_YZ - \nabla_Y\nabla_XZ - \nabla_{[X, Y]}Z, \eqno(A_1)$$
for all $X, Y, Z \in \Gamma(TM).$ First, by using (2.18a), (2.6b) and (2.18c) we deduce that

$$\begin{array}{lc}\nabla_{\frac{\delta }{\delta x^k}}\nabla_{\frac{\delta }{\delta x^j}}\frac{\delta }{\delta x^i}
= \nabla_{\frac{\delta }{\delta x^k}}\left(\bar{\Gamma}^{\;\;h}_{i\ k}{\frac{\delta }{\delta x^h}} + \Phi^{-2}
K_{lk}{\frac{\partial }{\partial x^0}}\right) = \frac{\delta\bar{\Gamma}^{\;\;h}_{i\ j}}{\delta x^k}
\frac{\delta }{\delta x^h}\vspace{2mm}\\ +
\bar{\Gamma}^{\;\;l}_{i\ j}\left(\bar{\Gamma}^{\;\;h}_{l\ k}\frac{\delta }{\delta x^h} + \Phi^{-2}
K_{lk}{\frac{\partial }{\partial x^0}}\right) - 2\Phi^{-2}K_{ij}c_k{\frac{\partial }{\partial x^0}} + \Phi^{-2}
\frac{\delta K_{ij}}{\delta x^k}{\frac{\partial }{\partial x^0}}\vspace{2mm}\\ + \Phi^{-2}K_{ij}\left(K_k^h
\frac{\delta }{\delta x^h}
 + c_k\frac{\partial }{\partial x^0}\right) = \left\{\frac{\delta\bar{\Gamma}^{\;\;h}_{i\ j}}{\delta x^k}
\frac{\delta }{\delta x^h} + \bar{\Gamma}^{\;\;l}_{i\ j}\bar{\Gamma}^{\;\;h}_{l\ k}\right.\vspace{2mm}\\ \left.
 + \Phi^{-2}
K_{ij}K_k^h\right\}\frac{\delta }{\delta x^h} + \Phi^{-2}\left\{\frac{\delta K_{ij}}{\delta x^k}
+ K_{lk}\bar{\Gamma}_{i\ j}^{\;\;l} - K_{ij}c_k\right\}\frac{\partial }{\partial x^0}.
\end{array}\eqno(A_2)$$
Then, by using (2.5a) and (2.18b), we obtain

$$\begin{array}{lc}\nabla_{[{\frac{\delta }{\delta x^k}}, {\frac{\delta }{\delta x^j}}]}\frac{\delta }{\delta x^i}
= 2\omega_{jk}K_i^h \frac{\delta }{\delta x^h} + 2b_i\omega_{jk}\frac{\partial }{\partial x^0}.
\end{array}\eqno(A_3)$$
Now, taking account of $(A_2)$ and $(A_3)$ in $(A_1)$, and using (3.3a) and the spatial covariant derivative
of the extrinsic tensor field, we infer that

$$\begin{array}{lc}R\left({\frac{\delta }{\delta x^k}}, {\frac{\delta }{\delta x^j}}, \frac{\delta }{\delta x^i}
\right) = \left\{\bar{R}^{\;\;h}_{i\ jk} + \Phi^{-2}\left(K_{ij}K_k^h - K_{ik}K_j^h\right)\right\}
\frac{\delta }{\delta x^h}\vspace{2mm}\\ + \left\{\Phi^{-2}\left(K_{ij\vert k} - K_{ik\vert j} + K_{ik}c_j -
K_{ij}c_k\right) - 2b_i\omega_{jk}\right\}\frac{\partial }{\partial x^0}.
\end{array}\eqno(A_4)$$
Similar calculations by using (2.18), (2.6d) and (2.5b) lead us to the following:

$$\begin{array}{lcr}\nabla_{\frac{\delta }{\delta x^k}}\nabla_{\frac{\partial }{\partial x^0}}\frac{\delta }
{\delta x^i}
= \left\{\frac{\delta K^h_i}{\delta x^k} + K_i^j\bar{\Gamma}_{j\ k}^{\;\;h} + b_iK_k^h\right\}\frac{\delta }
{\delta x^h}\vspace{2mm}\\ \hspace*{26mm} + \left\{\Phi^{-2}K_i^jK_{jk} + \frac{\delta b_i}{\delta x^k} + b_ic_k\right\}
\frac{\partial }{\partial x^0}.\end{array}\eqno(A_5)$$

$$\begin{array}{lc}\nabla_{\frac{\partial }{\partial x^0}}\nabla_{\frac{\delta }{\delta x^k}}\frac{\delta }{\delta x^i}
= \left\{\frac{\partial\bar{\Gamma}^{\;\;h}_{i\ k}}{\partial x^0} +
K_j^h\bar{\Gamma}^{\;\;j}_{i\ k} + K_{ik}b^h\right\}\frac{\delta }{\delta x^h}\vspace{2mm}\\ + \left\{\Phi^{-2}
\frac{\partial K_{ik}}{\partial x^0} + b_j\bar{\Gamma}^{\;\;j}_{i\ k} - \Psi\Phi^{-2}
K_{ik}\right\}\frac{\partial }{\partial x^0}.\end{array}\eqno(A_6)$$

$$\begin{array}{lc}\nabla_{[{\frac{\delta }{\delta x^k}}, {\frac{\partial }{\partial x^0}}]}\frac{\delta }{\delta x^i}
= - a_kK_i^h \frac{\delta }{\delta x^h} - b_ia_k\frac{\partial }{\partial x^0}.
\end{array}\eqno(A_7)$$
Then, by using $(A_5)$, $(A_6)$ and $(A_1)$, and taking into account (3.3b) and both covariant derivatives of the
extrinsic curvature tensor field, we deduce that

$$\begin{array}{lc}R\left({\frac{\delta }{\delta x^k}}, {\frac{\partial }{\partial x^0}}, \frac{\delta }{\delta x^i}
\right) = \left\{\bar{R}^{\;\;h}_{i\ 0k} + b_iK_k^h - b^hK_{ik}\right\}\frac{\delta }{\delta x^k} \vspace{2mm}\\
+ \left\{b_{i\vert k} + b_ib_k - \Phi^{-2}\left(K_{ik\vert 0} + K_{ij}K_k^j - \Psi K_{ik}\right)\right\}
\frac{\partial }{\partial x^0}.\end{array}\eqno(A_8)$$
Next, we consider the curvature tensor fields of type $(0,4)$ of the connections $\nabla$ and $\bar{\nabla}$,
denoted by $R$ and $\bar{R}$, and given by

$$\begin{array}{lc}(a) \ \ \ \ R(X, Y, Z, U) = g(R(X, Y, U), Z), \vspace{2mm}\\
(b) \ \ \ \ \bar{R}(X, Y, {\cal{S}}Z, {\cal{S}}U) = \bar{g}(\bar{R}(X, Y, {\cal{S}}U), {\cal{S}}Z),
\end{array}\eqno(A_9)$$
for all $X, Y, Z, U \in \Gamma(TM)$. Then we have the following components with respect to the threading frame:

$$\begin{array}{lc}(a) \ \ \ \ R_{iljk} =
R\left(\frac{\delta }{\delta x^k}, \frac{\delta }{\delta x^j}, \frac{\delta }{\delta x^l}, \frac{\delta }
{\delta x^i}\right) = \bar{g}_{lh}R_{i\ jk}^{\;h}, \vspace{2mm} \\
(b) \ \ \ \ R_{i0jk} =
R\left(\frac{\delta }{\delta x^k}, \frac{\delta }{\delta x^j}, \frac{\partial }{\partial x^0}, \frac{\delta }
{\delta x^i}\right) = - \Phi^{2}R_{i\ jk}^{\;0}, \vspace{2mm} \\
(c) \ \ \ \ R_{il0k} =
R\left(\frac{\delta }{\delta x^k}, \frac{\partial }{\partial x^0}, \frac{\delta }{\delta x^l}, \frac{\delta }
{\delta x^i}\right) = \bar{g}_{lh}R_{i\ 0k}^{\;h}, \vspace{2mm} \\
(d) \ \ \ \ R_{i00k} =
R\left(\frac{\delta }{\delta x^k}, \frac{\partial }{\partial x^0}, \frac{\partial }{\partial x^0},
\frac{\delta }{\delta x^i}\right) = -\Phi^{2}R_{i\ 0k}^{\;0}, \vspace{2mm} \\
(e) \ \ \ \ \bar{R}_{iljk} =
\bar{R}\left(\frac{\delta }{\delta x^k}, \frac{\delta }{\delta x^j}, \frac{\delta }{\delta x^l}, \frac{\delta }
{\delta x^i}\right) = \bar{g}_{lh}\bar{R}_{i\ jk}^{\;\;h}, \vspace{2mm} \\
(f) \ \ \ \ \bar{R}_{il0k} =
\bar{R}\left(\frac{\delta }{\delta x^k}, {\frac{\partial }{\partial x^0}}, \frac{\delta }{\delta x^l}, \frac{\delta }
{\delta x^i}\right) = \bar{g}_{lh}\bar{R}_{i\ 0k}^{\;\;h}.\end{array}\eqno(A_{10})$$
These local components are used in both Sections 4 and 5 in order to deduce the final form for the structure
equations and the local components of the Ricci tensor, respectively.

\section*{Appendix B}

We derive some useful formulas for the conservation equations stated in Section 8. First, by using (6.10), (2.18d)
 and (6.13a), we obtain

$$\begin{array}{lc} (\nabla_UT)(U, \frac{\partial }{\partial x^0}) = (\nabla_{\Phi^{-1}\frac{\partial }{\partial x^0}}T)
(\Phi^{-1}\frac{\partial }{\partial x^0}, \frac{\partial }{\partial x^0})\vspace{2mm}\\
 = \Phi^{-2}\left\{\frac{\partial T_{00}}{\partial x^0} - 2T(\nabla_{\frac{\partial }{\partial x^0}}\frac{\partial }{\partial x^0},
\frac{\partial }{\partial x^0})\right\}\vspace{2mm}\\  = \Phi^{-2}\left\{\frac{\partial T_{00}}{\partial x^0}
- 2T(\Phi^2b^k\frac{\delta }{\delta x^k} + \Psi\frac{\partial }{\partial x^0}, \frac{\partial }
{\partial x^0})\right\}\vspace{2mm}\\
 = \Phi^{-2}\left\{\frac{\partial T_{00}}{\partial x^0} - 2\Psi T_{00}\right\} - 2b^jT_{j0}
.\end{array}\eqno(B_1)$$
Then, by using (5.2), (6.13a), (5.3), (2.18a), (6.10) and (6.15a), we deduce that

$$\begin{array}{r}\displaystyle{\sum_{k=1}^3}\left\{(\nabla_{E_k}T)(E_k, \frac{\partial }{\partial x^0})\right\} =
\displaystyle{\sum_{k=1}^3}\left\{(\nabla_{E_k^h\frac{\delta }{\delta x^h}}T)(E_k^j\frac{\delta }{\delta x^j},
 \frac{\partial }{\partial x^0})\right\}\vspace{2mm}\\ =
\displaystyle{\sum_{k=1}^3}\left\{E_k^h\frac{\delta }{\delta x^h}(E_k^jT_{j0})
- E_k^hT(\nabla_{\frac{\delta }{\delta x^h}}(E_k^j\frac{\delta }{\delta x^j}), \frac{\partial }{\partial x^0})
\right.\vspace{2mm}\\ \left.- E_k^hE_k^jT(\frac{\delta }{\delta x^j}, \nabla_{\frac{\delta }{\delta
 x^h}}\frac{\partial }{\partial  x^0})\right\} = \bar{g}^{jh}\left\{ \frac{\delta T_{j0}}{\delta x^h}
- T(\bar{\Gamma}_{j\ h}^{\;k} \frac{\delta }{\delta x^k} \right.\vspace{2mm}\\ \left. + (\omega_{jh} + \Phi^{-2}
 \Theta_{jh})\frac{\partial }{\partial x^0},  \frac{\partial }{\partial x^0})
- T(\frac{\delta }{\delta x^j},
 (\Theta_h^k + \Phi^2\omega_h^k)\frac{\delta }{\delta x^k} + c_k\frac{\partial }{\partial x^0})\right\}\vspace{2mm}\\
 = T^j_{0\vert j} - \Phi^{-2}\Theta T_{00} - \Theta^{jk}T_{jk} - c^jT_{j0}
.\end{array}\eqno(B_{2})$$
Next, by similar calculations, we infer that

$$\begin{array}{lc} (\nabla_UT)(U, \frac{\delta }{\delta x^i}) = \Phi^{-2}\left\{T_{i0\vert 0} -
\Psi T_{i0} - b_iT_{00}\right\} - T_{ij}b^j,
\end{array}\eqno(B_{3})$$
and

$$\begin{array}{r}\displaystyle{\sum_{k=1}^3}\left\{(\nabla_{E_k}T)(E_k, \frac{\delta }{\delta x^i})\right\}
= T^j_{i\vert j} - \Phi^{-2}\Theta T_{i0} - (\omega_{ij} + \Phi^{-2}\Theta_{ij})T^j_0
.\end{array}\eqno(B_4)$$

\section*{Appendix C}

In this appendix we present the calculations for both the spatial Ricci tensor and the spatial scalar curvature of
an almost FLRW universe. First, by using (2.16), (9.4a), (9.4d) and (9.4e), we obtain

$$\bar{\Gamma}^{\;\;k}_{i\ j} = \frac{1}{1 - 2B}\left\{\delta_{ij}B^k - \delta_i^kB_j - \delta_j^kB_i\right\},
 \eqno(C_1)$$
 where we put $B^k = B_k = \partial B/\partial x^k$. Then (2.15a) and $(C_1)$ imply

$$\begin{array}{lc}\bar{g}\left(\bar{\nabla}_{\frac{\partial }{\partial x^j}}\frac{\partial }{\partial  x^i},
 \frac{\partial }{\partial x^h}\right)
= a^2\left\{\delta_{ij}B_h - \delta_{ih}B_j - \delta_{jh}B_i\right\}
.\end{array}\eqno(C_2)$$
Then, apply $\partial/\partial x^k$ to both sides in $(C_2)$ and taking into account that $\bar{\nabla}$ is a
metric connection, we deduce that

$$\begin{array}{lc}\bar{g}\left(\bar{\nabla}_{\frac{\partial }{\partial x^k}}\bar{\nabla}_{\frac{\partial }
{\partial x^j}}\frac{\partial }{\partial  x^i}, \frac{\partial }{\partial x^h}\right) +
\bar{g}\left(\bar{\nabla}_{\frac{\partial }{\partial x^j}}\frac{\partial }{\partial  x^i},
\bar{\nabla}_{\frac{\partial }{\partial x^k}}\frac{\partial }{\partial x^h}\right)\vspace{2mm}\\
= a^2\left\{\delta_{ij}B_{hk} - \delta_{ih}B_{jk} - \delta_{jh}B_{ik}\right\},\end{array}\eqno(C_3)$$
where we put

$$B_{jk} = \frac{\partial^2B}{\partial x^j\partial x^k}.$$
Denote the left hand side of $(C_3)$ by $L_{ihjk}$ and by using (3.1), (3.2a), (2.15a) and $(C_1)$ we infer that

$$\begin{array}{r}L_{ihjk} - L_{ihkj} = R_{i\ jk}^{\;l}\bar{g}_{lh} + \left(\bar{\Gamma}_{i\ j}^{\;l}
\bar{\Gamma}_{h\ k}^{\;m} - \bar{\Gamma}_{i\ k}^{\;l}\bar{\Gamma}_{k\ j}^{\;m}\right)\bar{g}_{lm}
= \bar{R}_{ihjk}\vspace{2mm}\\ +  \frac{a^2}{1 - 2B}\left\{\left(\delta_{ij}B^l - \delta_i^lB_j
- \delta_j^lB_i\right)\left(\delta_{hk}B_l- \delta_{hl}B_k - \delta_{kl}B_h\right)\right.\vspace{2mm}\\ \left.
 - \left(\delta_{ik}B^l - \delta_i^lB_k - \delta_k^lB_i\right)\left(\delta_{hj}B_l - \delta_{hl}B_j
- \delta_{jl}B_h\right)\right\}.\end{array}\eqno(C_4)$$
Thus, $(C_3)$ and $(C_4)$ imply

$$\begin{array}{r}\bar{R}_{ihjk} = a^2\left\{\delta_{ij}B_{hk} + \delta_{hk}B_{ij} - \delta_{ik}B_{hj} -
 \delta_{jh}B_{ik} \right\}\vspace{2mm}\\ +  \frac{a^2}{1 - 2B}\left\{\left(\delta_{ik}b^l - \delta_i^lB_k
- \delta_k^lB_i\right)\left(\delta_{hj}B_l- \delta_{hl}B_j - \delta_{jl}B_h\right)\right.\vspace{2mm}\\ \left.
 - \left(\delta_{ij}b^l - \delta_i^lB_j - \delta_j^lB_i\right)\left(\delta_{hk}B_l - \delta_{hl}B_k
- \delta_{kl}B_h\right)\right\}.\end{array}\eqno(C_5)$$
Now, contracting $(C_5)$ by $\bar{g}^{hk}$ and using (5.9), (5.12) and (9.4e), we obtain

$$\begin{array}{r}\bar{R}_{ij} = \frac{1}{(1 - 2B)^2}\left\{\delta_{ij}\displaystyle{\sum_{k=1}^3}(B_k)^2 +
3B_iB_j + (1 - 2B)\left(B_{ij} + \delta_{ij}\displaystyle{\sum_{k=1}^3}B_{kk}\right)\right\}
,\end{array}\eqno(C_6)$$
and
$$\bar{\bf R} = \frac{2}{a^2(1 - 2B)^3}\left\{3\displaystyle{\sum_{k=1}^3}(B_k)^2 +
2(1 - 2B)\displaystyle{\sum_{k=1}^3}B_{kk}\right\}.\eqno(C_7)$$

\end{document}